\documentclass[10pt,a4paper]{article}
\usepackage{amsmath}
\usepackage{amssymb}
\usepackage{theorem}
\theoremheaderfont{\scshape}
\newtheorem{thm}{Theorem}[section]
\newtheorem{cor}[thm]{Corollary}
\newtheorem{lem}[thm]{Lemma}

\theorembodyfont{\normalfont}
\newtheorem{rem}[thm]{Remark}

\def\bb#1{\mathbb{#1}}
\def\cal#1{\mathcal{#1}}
\def\ds{\displaystyle}

\def\eps{\varepsilon}

\def\imp{\Rightarrow}
\def\ben{\begin{enumerate}}
\def\een{\end{enumerate}}

\def\limf#1#2{\ensuremath{\underset{#1\rightarrow #2}{\lim} }}

\def\abs#1{\left\vert #1\right\vert}
\def\N#1{\|#1\|}

\def\interventier#1#2{\ensuremath{[\![#1,#2]\!]}}

\newcommand{\Sum}[2]{\displaystyle \sum\limits_{#1}^{#2}}

\newcommand{\Lim}[2]{\lim\limits_{#1\rightarrow #2}}

\title{Exponential inequalities for martingales and \\
asymptotic properties  of the free energy of \\ directed polymers in
random environment}
\author{Quansheng Liu\thanks {
Email: Quansheng.Liu@univ-ubs.fr}
              \;\;and \;\;  Fr\'ed\'erique
              Watbled\thanks{Email: Frederique.Watbled@univ-ubs.fr}\\
\\
LMAM, Universit\'e de Bretagne Sud,\\
Campus de Tohannic,
BP 573, 56017 Vannes, France\\
Universit\'e Europ\'eenne de Bretagne}

\date{November 2008}

\begin{document}

\maketitle


\begin{abstract}
We first obtain exponential inequalities for martingales. Let $(X_k)
(1\leq k\leq n)$ be a sequence of martingale differences relative to
a filtration $(\cal F_k)$, and set $S_n=X_1 + ... + X_n$.
 We prove that if for some $\delta>0, Q\geq 1$, $K>0$ and all $k$, a.s.
 $\bb E[e^{\delta\abs{X_k}^Q}|\cal F_{k-1}]\leq K$, then for some constant
$c>0$ (depending only on $\delta, Q$ and $K$) and all $x>0$,
 $P[|S_n|>nx]\leq 2 e^{-nc(x)}$, where $c(x) = cx^2$ if $x\in
 ]0,1]$, and $c(x) = cx^Q$ if $x>1$; the converse also holds if
 $(X_i)$ are independent and identically distributed. It extends
  Bernstein's inequality for $Q=1$, and Hoeffding's inequality for
$Q=2$. We then apply the preceding result to establish exponential
concentration inequalities for the free energy of directed polymers
in random environment, show its rate of convergence (in probability,
almost surely, and in $L^p$), and
give it an expression in terms of free energies of some
multiplicative cascades, which improves an inequality of Comets and
Vargas (2006, \cite{CometsVargas}) to an equality.
\end{abstract}

\medskip
{\bf Key words}. Martingale differences,  super-martingales, large
deviation inequality, exponential inequality, Bernstein's inequality,
Hoeffding-Azuma's inequality, directed polymers, random environment,
concentration inequality, free energy, convergence rate,
multiplicative cascades.

\medskip
{\bf 2000 AMS Subject Classifications.} Primary 60G42,60K35;
secondary 60K37, 60G50, 60F10, 82D30

\bigskip

\section{Introduction and main results}
Our work was initially motivated by the study of the free energy of
a directed polymer in a random environment. Comets and Vargas (2006, \cite{CometsVargas})
proved that the free energy (at $\infty$) is bounded by the infimum
of those of some generalized multiplicative cascades, and that the
equality holds if the environment is bounded or gaussian. The
essential point in their proof for the equality is an exponential
concentration inequality for the free energy (at time $n$), which
was not known for a general environment. Using a large deviation
inequality of Lesigne and Volny (2001, \cite{LesigneVolny}) on martingales, Comets, Shiga
and Yoshida (2003, \cite{CometsShigaYoshida}) did obtain a concentration inequality for the
free energy;  but their bound is larger than the exponential one,
and is not sharper enough to imply the equality mentioned above.
Another non satisfactory point of their inequality is that it cannot be used to prove
rigourous results on the rate of convergence, for the almost sure (a.s.)
or $L^p$ convergence of the free energies.

\medskip
The objective of the present paper is to establish exponential large
deviation inequalities, and to use
 them  to show exponential concentration inequalities for the free energy of a
polymer in general random environment, its rate of convergence, and
an expression of its limit value in terms of those of some
multiplicative cascades.

\medskip
Large deviation inequalities are very powerful tools in probability
theory, and have been studied by many authors: see e.g. the
classical works of Bernstein (1924, \cite{Bernstein1924}), Cram\'er (1938, \cite{Cramer}), Hoeffding
(1963, \cite{Hoeffding}), Azuma (1967, \cite{Azuma}), Chernoff (1981, \cite{Chernoff}),
the books of Chow and
Teicher (1978, \cite{ChowTeicher}), and  Petrov (1995, \cite{Petrov}), and the recent papers by de la
Pe\~na, (1999, \cite{delaPena}), Lesigne and Voln\'y (2001, \cite{LesigneVolny}),
Bentkus (2004, \cite{Bentkus2004}),
and Chung and Lu (2006, \cite{ChungLu2006}).
See also Ledoux  (1999, \cite{Ledoux1999}) and Wang (2005, \cite{Wang2005}) for related concentration inequalities and general functional inequalities.

 Let $(\Omega,\cal F,P)$ be a probability
space, and let $\cal F_0=\{\emptyset,\Omega\}\subset \cal
F_1\subset\cdots\subset \cal F_n$ be an increasing sequence of
sub-$\sigma$-fields of $\cal F$. Let $ X_1, ..., X_n$ be a sequence
of real- valued martingale differences defined on $(\Omega,\cal
F,P)$, adapted to the filtration $(\cal F_k)$: that is, for each
$1\leq k \leq n$,   $X_k$ is $\cal F_k$ measurable, and $\bb
E[X_k|\cal F_{k-1}] = 0$. Set
\begin{equation}
 S_n = X_1 + ... + X_n.
 \end{equation}
We are interested in exponential large deviation inequalities of the form
\begin{equation}
\label{cramer-bound}
 P[|S_n| > nx]= O (e^{-c(x)n}),
 \end{equation}
where $x>0$ and $c(x)>0$.  When $(X_i)$ are independent and
identically distributed (iid) with mean $\bb EX_i=0$, it is known [see Petrov
(1995, \cite{Petrov} p.137)] that (\ref{cramer-bound}) holds for all $x>0$ and
some $c(x)>0$ if and only if for some $\delta >0$,
\begin{equation}
\label{cramer-cond}
  \bb E e^{\delta |X_1| } <\infty.
\end{equation}
For a sequence of  martingale differences, Lesigne and Voln\'y
(2001, \cite{LesigneVolny}) proved that if for some constant $K>0$ and all $k=1, ..., n$,
 \begin{equation}
 \label{exp-moment}
       \bb E e^{ |X_k| } \leq K,
\end{equation}
then for any $x>0$,
\begin{equation}
\label{LV}
 P\left[\frac{S_n}{n} > x\right] = O ( e^{ - \frac{1}{4}  x^{2/3} n^{1/3}}),
\end{equation}
and that this is the best possible inequality that we can have under
the condition (\ref{exp-moment}),  even in the class of stationary and ergodic sequences of martingale differences,
in the sense that there exist such
sequences of martingale differences $(X_i)$
satisfying (\ref{exp-moment}) for some $K>0$, but
\begin{equation}
\label{LV2}
 P\left[\frac{S_n}{n} > 1\right]> e^{ - c n^{1/3}}
 \end{equation}
for some constant $c>0$ and infinitely many $n$. It is therefore interesting to know what is
the good condition to have the exponential inequality
(\ref{cramer-bound}) in the martingale case. It turns out that
(\ref{cramer-bound}) still holds if we replace the expectation in
(\ref{exp-moment}) by the conditional one given $\cal
F_{k-1}$. In fact we shall prove the following much sharper result.
It is a consequence of Theorems \ref{superdiff1}, \ref{suiteadaptee}, and \ref{superdiffQ}.

\begin{thm}\label{supermartdiff}
Let $(X_k)$ be a $\{ \cal F_k \}$-adapted sequence of martingale
differences.   Assume that for some constants $Q\geq 1$, $\delta>0$,
$K>0$ and all $k\in\{1,\cdots,n\}$, almost surely
\begin{equation}\label{C}
\bb E[e^{\delta\abs{X_k}^Q}|\cal F_{k-1}]\leq K.
\end{equation}
Then there exists a constant $c>0$ depending only on $Q$, $\delta$
and $K$, such that for all $x>0$,
\begin{equation}\label{th-intro1}
P\left[\pm\frac{S_n}{n}>x \right]\leq
 \left\{\begin{aligned}
& e^{-ncx^2}  \;\; \textrm{ if }  x\in ]0,1], \\
& e^{-ncx^Q}   \;\; \textrm{ if } x\in ]1,\infty[.
\end{aligned}\right.
\end{equation}
The converse also holds in the iid case: if $X_k$ are iid and if
$P[\pm\frac{S_n}{n}>x ]\leq e^{-ncx^Q}$ holds for some $n\geq 1$,
$Q\geq 1$, $c>0$, $x_1>0$ and all $x\geq x_1$, then for all $\delta
\in ]0,c[$, there exists $K=K(\delta,Q,c,x_1)>0$ such that
$$\bb E[e^{\delta\abs{X_1}^Q}]\leq K.$$
\end{thm}

By the result of  Lesigne and Voln\'y
(\cite{LesigneVolny}) cited above, the conditional exponential moment
condition  (\ref{C}) cannot be relaxed to the non conditional one.

When $(X_k)$ are iid with $\bb E[X_k]=0$, Bernstein's inequality states
(cf. \cite{Petrov}, page 57) that
if $\sigma^2=\bb E[X_k^2]<\infty$ and
\begin{equation}\label{BC}
\abs{\bb EX_k^m}\leq \frac{1}{2}m!\sigma^ 2H^ {m-2}
\end{equation}
for some $H>0$ and all $m=2,3,\cdots$, then
\begin{equation}\label{BI}
P\left[\pm\frac{S_n}{n}>x \right]\leq
 \left\{\begin{aligned}
& e^{-nc_0x^2}  \;\; \textrm{ if }  x\in ]0,x_0], \\
& e^{-nc_1x}   \;\; \textrm{ if } x\in ]x_0,\infty[,
\end{aligned}\right.
\end{equation}
where $c_0=\frac{1}{4\sigma^ 2}$, $c_1=\frac{1}{4H}$,
$x_0=\frac{\sigma^ 2}{H}$.
Notice that (in the iid case) Bernstein's condition (\ref{BC}) is
equivalent to Cramer's condition that $\exists \delta >0$ such that
\begin{equation}\label{cramer}
\bb E[e^{\delta\abs{X_k}}]<\infty.
\end{equation}
In applications we find more convenient to use Cramer's condition.
Taking $Q=1$ in Theorem \ref{supermartdiff}, we obtain the following
Bernstein-type inequality.

\begin{cor} (A Bernstein-type inequality)
Assume that $(X_k)$ are iid with $\bb E[X_k]=0$, $1\leq k\leq n$.
If (\ref{cramer}) holds for some $\delta>0$, then for some $c=c(\delta)>0$,
\begin{eqnarray}\label{}
P\left[\pm\frac{S_n}{n}>x \right]\leq
 \left\{\begin{aligned}
& e^{-ncx^2}  \;\; \textrm{ if }  x\in ]0,1], \\
& e^{-ncx}   \;\; \textrm{ if } x\in ]1,\infty[.
\end{aligned}\right.
\end{eqnarray}
Conversely, if for some $n\geq 1$, $c>0$,  $x_0>0$ and all $x>x_0$,
$
  P\left[\pm\frac{S_n}{n}>x \right]\leq e^{-ncx},    $
 then (\ref{cramer}) holds for each $\delta \in
]0,c[$.
\end{cor}

When $Q=2$, Theorem \ref{supermartdiff} extends the following
well-known Hoeffding's inequality{\footnote {The inequality
(\ref{Hoeffding}) is often called Hoeffding's inequality when
$(X_k)$  are iid, and Azuma's inequality when $(X_k)$ are martingale
differences. This is rather strange, as it was Hoeffding (1963) who
first obtained it for martingales, although he mainly treated the
iid case, and only mentioned the martingale case as a remark [see
\cite{Hoeffding}, p.18]. To respect the history, we call it
Hoeffding's inequality, although Azuma (1967, \cite{Azuma}) refound
it four years later. We think that what happened would be that, the
first author who called it Azuma's inequality did not know the
existence of the  remark of Hoeffding, the second followed the first
without verification, and so on. } }:
 if $(X_k)$ is a sequence of martingale
differences with $|X_k| \leq a$ a.s. for some constant $a\in
]0,\infty[$, then for all $n \geq  1$ and all $x>0$,
\begin{equation}\label{Hoeffding}
 P\left[\pm\frac{S_n}{n} > x \right]\leq e^{-ncx^2},
\end{equation}
where $c=1/(2a^2)$. In fact, by our result for $Q=2$, we obtain:

\begin{cor} \label{supermartdiffcor}(Extension of Hoeffding's inequality)
When $(X_k)$ are iid,  then there is a constant   $c>0$ such that
(\ref{Hoeffding}) holds for all $n\geq 1$ and all $x>0$,   if and
only if for some $\delta >0$,
\begin{equation}
 \bb Ee^{\delta X_1^2} < \infty.
\end{equation}
Moreover, if (\ref{Hoeffding}) holds for some $n\geq 1$ and all $x>0$,
with some constant $c=c_1$, then it holds for all $n \geq 1$ and $x>0$,
with some constant $c=c_2$ depending only on $c_1$.
\end{cor}

So our result is a complete extension of Hoeffding's inequality
even in the iid case.

\medskip
We then apply the preceding results to directed polymers in random
environment that we describe as follows.
 Let $(\omega_n)_{n\in\bb
N}$ be the simple random walk on $\bb Z^d$ starting at $0$, defined
on a probability space $(\Omega,\cal F,P)$. Let
$(\eta(n,x))_{(n,x)\in\bb N\times \bb Z^d}$ be a sequence of i.i.d.
real random variables defined on another probability space $(E,\cal
E, \bb Q)$ (we use the letter E to refer the Environment).  For real
$\beta$ (the inverse of temperature), define
\begin{equation}
\lambda(\beta)=\ln \bb Q[e^{\beta\eta(0,0)}].
\end{equation}
(If $\mu$  is a measure and f is a function, we write $\mu(f)$ or   $\mu[f]$ for the
integral of $f$ with respect to $\mu $.)
 We  fix  $\beta >0$, and only suppose that
 \begin{equation}\label{H}
 \bb Q [e^{\beta|\eta(0,0)|}]<\infty
 \end{equation}
 (we do not suppose that it holds for all $\beta>0$).
 Of course this condition is equivalent to $\lambda (\pm \beta) < \infty$.
 We are interested in the normalized partition function
\begin{equation}
W_n(\beta)=P\left[\exp\left(\beta\sum_{j=1}^n \eta(j,\omega_j) -n\lambda(\beta)\right)\right],
\end{equation}
and the free energy $\frac{1}{n} \ln W_n(\beta)$.

This model first appeared in physics literature [see Huse and Henley
(1985, \cite{HuseHenley})]  to modelize the phase boundary of Ising model subject to
random impurities; the first mathematical study was undertaken by
Imbrie and Spencer (1988, \cite{ImbrieSpencer}) and Bolthausen (1989, \cite{Bolthausen}). For recent results,
see e.g. Carmona and Hu (2004, \cite{CarmonaHu2004}),
Carmona, Guerra, Hu and M\'ejane (2006, \cite{CarmonaGuerraHuMejane2006}), Comets, Shiga and Yoshida (2004, \cite{CometsShigaYoshidaReview}),
and Comets and Yoshida (2006, \cite{CometsYoshida2006}).

\medskip
Assuming $Q [e^{\beta|\eta(0,0)|}] <\infty$ {\em for all } $\beta
>0$,
 Comets, Shiga and Yoshida (\cite{CometsShigaYoshida}) proved that $\forall x>0$, there
exists $n_0\in\bb N^*$ such that for any $n\geq n_0$,
\begin{equation} \label{cometsconcentration}
\bb Q\left[\abs{\frac{1}{n}\ln W_n(\beta)-\frac{1}{n}\bb Q[\ln
W_n(\beta)]}>x\right]\leq
\exp\left(-\ds\frac{n^{\frac{1}{3}}x^{\frac{2}{3}}}{4}\right).
\end{equation}
In fact, in their proof of (\ref{cometsconcentration}), they used
the condition that $Q [e^{3\beta|\eta(0,0)|}] <\infty$, due to the
application of their Lemma 3.1 (p.711).

\medskip
We first improve this result to an exponential inequality under the
weaker condition that $ Q [e^{\beta|\eta(0,0)|}] < \infty$ for the
fixed $\beta$.

\begin{thm}\label{polymeres} (Exponential concentration inequality for the free
energy) Let $\beta>0$ be fixed such that $  Q [e^{\beta|\eta(0,0)|}]
< \infty$.  If for some $Q\geq 1$ and $R>0$,
\begin{equation}
\label{condQ}
 \bb Q[e^{R\abs{\eta(0,0)}^Q}]<+\infty,
\end{equation}
then
\begin{equation}\label{tailQ}
\bb Q\left[\frac{1}{n}\abs{\ln W_n(\beta)-\bb Q[\ln W_n(\beta)]}>x\right]\leq \left\{\begin{aligned}
&2e^{-ncx^2} \textrm{ if } 0\leq x\leq 1,\\
&2e^{-n cx^Q} \textrm{ if } x>1,\\
\end{aligned}\right.
\end{equation}
where $c>0$ is a constant depending only on $Q, R$, and the law of $\eta(0,0)$.
\end{thm}
Notice that the condition (\ref{condQ}) holds automatically for
$Q=1$ and $R=\beta$, so that $(\ref{tailQ})$ holds for $Q=1$ under
the only hypothesis $ Q [e^{\beta|\eta(0,0)|}] < \infty$; when
(\ref{condQ}) holds  for some $Q>1$ and $R>0$,  $(\ref{tailQ})$
gives a sharper bound for large values of $x$.

Theorem  \ref{polymeres} is a consequence of Corollary
\ref{corollaireQ}. As shown in Carmona and Hu (2002,
\cite{CarmonaHu2002}) and Comets and Vargas (2006,
\cite{CometsVargas}), when the environment is gaussian or bounded,
the inequality can be obtained directly by a general concentration
result on gaussian or bounded variables (see e.g. Ledoux (1999,
\cite{Ledoux1999})). But this method does not work for a general
environment.

As applications we shall show the following properties about the
free energy $\frac{1}{n}\ln W_n(\beta)$:
\begin{list}{}
\item (1) for some $p_{-}(\beta)\leq 0$, $\frac{1}{n}\ln W_n(\beta)\to p_{-}(\beta)$ in probability \textit{at an exponential rate} (cf. Theorem \ref{energy});
\item (2) $\frac{1}{n}\ln W_n(\beta)\to p_{-}(\beta)$ a.s. and in $L^p$, for all $p\geq 1$, at a rate
$O\left(\sqrt{\frac{\ln n}{n}}\right)$ (cf. Theorem \ref{estimationdetermin});
\item (3) $p_{-}(\beta)$ can be expressed in terms of some generalized multiplicative cascades
(cf. Theorem \ref{cascades}).
\end{list}
\medskip Part (1) extends the same conclusion of Carmona and Hu
(2002, \cite{CarmonaHu2002}) for the gaussian environment case to a
general environment case. The rate of a.s. convergence in part (2)
improves the bound $O(n^{-(\frac{1}{2}-\eps)})$ ($\eps>0$) of
Carmona and Hu (2004, \cite{CarmonaHu2004}) obtained for the
gaussian environment case. Part (3) improves an inequality of Comets
and Vargas (2006, \cite{CometsVargas}) to an equality.

\medskip
The rest of the paper is organized as follows. In Section $2$ we
establish exponential inequalities for supermartingales, which
extend Bernstein or  Hoeffding's inequalities, according to
 $\bb E\left[e^{\delta\abs{X_i}}|\cal F_{i-1}\right]\leq K$ or
   $\bb E\left[e^{\delta\abs{X_i}^2}|\cal F_{i-1}\right]\leq K$,
   respectively. For large values of $x$, sharper inequalities are proven in
   Section $3$ under the condition that $\bb E\left[e^{\delta\abs{X_i}^Q}|\cal F_{i-1}\right]\leq
   K$ ($Q>1$).   These results are extended in Section $4$ to the more general
   case where $\bb E\left[e^{\delta\abs{X_i}^Q}|\cal F_{i-1}\right]\leq
   K_i$. As applications,  we show in Section $5$ the rate of
   convergences, a.s. and in $L^p$.
In the last $3$ sections, we study the free energies of directed
polymers in
   random environment, with the help of
our results on martingales:  we show exponential concentration
inequalities for the free energies in Section
   $6$, their convergence rates (in probability, a.s. and in $L^p$)
   in Section $7$, and, in Section 8, an expression of their limit value in terms of
   some generalized multiplicative cascades.

\section{Exponential inequalities for supermartingales}
\setcounter{equation}{0}

In this section we give an extension of Bernstein and Hoeffding's
inequalities to supermartingales with unbounded differences. Our
results are sharp even in the iid case.

Let $(X_i)_{1\leq i \leq n}$ be a sequence of real-valued
supermartingale differences defined on a probability space
$(\Omega,\cal F,P)$, adapted to a filtration $(\cal F_i)$, with
$\cal F_0=\{\emptyset,\Omega\}$. This means that for each
$1\leq i \leq n$, $X_i$ is $\cal F_i$-measurable, and
$\bb E[X_i|\cal F_{i-1}]\leq 0$ a.s..
We are interested in the growth rate of the Laplace transform $\bb
E[e^{tS_n}]$, and the convergence rate of the deviation
probabilities $P\left[\frac{S_n}{n}>x\right]$.

   \begin{thm} \label{superdiff1}
Let $(X_i)_{1\leq i\leq n}$ be a finite sequence of supermartingale
differences. If for some constant $K>0$ and all $i\in \interventier{1}{n}$,
\begin{equation}\label{K}
\bb E[e^{\abs{X_i}}|\cal F_{i-1}]\leq K\ \ \ a.s.,
\end{equation}
then:
\begin{equation}\label{growthsuperdiff1}
 \bb E[e^{tS_n}]\leq \exp\left({\frac{nKt^2}{1-t}}\right)
\textrm{ for all }t\in ]0,1[,
\end{equation}
and \begin{equation}\label{concsuperdiff1}
P\left[\frac{S_n}{n}>x\right]\leq\exp\left(-n\left(\sqrt{x+K}-\sqrt{K}\right)^2\right) \textrm{ for all } x>0.
\end{equation}
Consequently, \begin{equation}\label{concsuperdiff1bis}
 P\left[\frac{S_n}{n}>x\right]\leq \left\{\begin{aligned}
&\exp\left(-\ds\frac{nx^2}{K(1+\sqrt{2})^2}\right) &\textrm{ if } x\in ]0,K],\\
&\exp\left(-\ds\frac{nx}{(1+\sqrt{2})^2}\right)  &\textrm{ if } x\in
]K, \infty[.
\end{aligned}\right.
\end{equation}
Conversely, if $(X_k)$ are iid, and if
$P\left[\frac{S_n}{n}>x\right]\leq e^{-ncx}$ for some $n>1$, $c>0$,
and all $x\geq x_1>0$ large enough, then for all $\delta \in ]0,c[$,
$$\bb E[e^{\delta {X_1^+}}]\leq K,
  \;\;  \mbox{ where }  X_1^+ = \max (X_1, 0), \textrm{ and } K=e^{\delta x_1}+\frac{\delta}{c-\delta}e^{-(c-\delta)x_1}.$$
\end{thm}

\begin{cor}\label{cor1superdiff1}
Under the conditions of Theorem \ref{superdiff1}, $\forall\eps >0$, there exist $0<x_0<x_1$ and $K_1>0$ depending only on $K$ and $\eps$, such that:
\begin{equation}\label{corconcsuperdiff1}
P\left[\frac{S_n}{n}>x\right] \leq
\left\{\begin{aligned}
&\exp\left(-\ds\frac{nx^2}{4K(1+\eps)}\right)\textrm{ if } x\in ]0,x_0[,\\
&\exp\left(-\ds\frac{nx}{K_1}\right)\textrm{ if } x\in [x_0,x_1],\\
&\exp\left(-\ds\frac{nx}{1+\eps}\right)\textrm{ if } x\in ]x_1,+\infty[.
\end{aligned}\right.
\end{equation}
\end{cor}

We divide the proof into a series of lemmas.

   \begin{lem}\label{lemsuperdiff1.1}
Let $(X_i)_{1\leq i\leq n}$ be a finite sequence of random variables
adapted to a filtration $(\cal F_i)_{1\leq i\leq n}$. Let
$(l_i)_{1\leq i\leq n}$ be a finite sequence of deterministic
functions defined on a subinterval $I$ of $]0,\infty[$, such that
for each $i$ and each $t \in I$,
\begin{equation}\label{cond-lem1.1}
\bb E[e^{tX_i}|\cal F_{i-1}]\leq e^{l_i(t)}\ \ \  a.s..
 \end{equation}
 Then for every  $t\in I$,
\begin{equation}\label{cond-lem1.2}
\bb E[e^{tS_n}]\leq \exp\left(\sum_{i=1}^n l_i(t)\right),
 \end{equation}
and for  every $x>0$,
\begin{equation}\label{cond-lem1.3}
P\left[\frac{S_n}{n}>x\right]\leq e^{-nL_n^*(x)},
\end{equation}
where
\begin{equation}\label{cond-lem1.4}
L_n(t)=\frac{1}{n}\sum_{i=1}^n l_i(t), \;\; \mbox{ and }
   L_n^*(x)=\sup_{t\in I}\left(tx-L_n(t)\right).
\end{equation}
   \end{lem}

\textbf{Proof.}
(\ref{cond-lem1.2}) can be obtained by a simple induction argument on $n$.
(\ref{cond-lem1.3}) is an immediate consequence of (\ref{cond-lem1.2}), since $\forall x>0$,
$\forall t\in I$,
$$P\left[\frac{S_n}{n}>x\right]=P[e^{tS_n}> e^{tnx}]
\leq e^{-ntx}\bb E[e^{tS_n}]\leq \exp\left(-n(tx-L_n(t))\right).$$
\hfill{\rule{2mm}{2mm}\vskip3mm \par}

\begin{rem}
The submultiplicativity (\ref{cond-lem1.2}) for an adapted sequence corresponds to the multiplicativity
$\bb E[e^{tS_n}]=\prod_{i=1}^n \bb E[e^{tX_i}]$ in the independent case.
This explains why it is natural to consider the conditional Laplace transform
$\bb E[e^{tX_i}|\cal F_{i-1}]$ in the supermartingale case, instead of the Laplace transform
$\bb E[e^{tX_i}]$ in the independent case. For example, using Lemma \ref{lemsuperdiff1.1},
we can obtain the following generalization of Petrov's inequality (p.54 of \cite{Petrov}):

\begin{lem}\label{lempetrov}
Let $a_i>0$ and $T>0$ be constants such that for all $1\leq i\leq n$ and all $t\in ]0,T]$,
a.s. $\bb E[e^{tX_i}|\cal F_{i-1}]\leq e^{a_it^2}$.
Then for each $A\geq\frac{1}{n}\Sum{i=1}{n}a_i$, we have
\begin{equation}\label{concpetrov}
P\left[\frac{S_n}{n}>x\right] \leq
\left\{\begin{aligned}
&\exp\left(-\ds\frac{nx^2}{4A}\right)\textrm{ if } x\in ]0,2AT[,\\
&\exp\left(-\ds\frac{nTx}{2}\right)\textrm{ if } x\in [2AT,+\infty[.
\end{aligned}\right.
\end{equation}
\end{lem}

\textbf{Proof.} We apply Lemma \ref{lemsuperdiff1.1} with $I=]0,T]$ and
$l_i(t)=a_it^2$, $L_n(t)=At^2$,
which gives:
$$\bb E[e^{tS_n}]\leq \exp(nAt^2) \textrm{ for every } t\in ]0,T],$$
and
$$P\left[\frac{S_n}{n}>x\right]\leq e^{-nL_n^*(x)},$$
with $L_n^*(x)=\sup_{t\in ]0,T]}\left(tx-At^2\right)$.
We calculate this $\sup$ and find:
$$L_n^*(x)=\left\{\begin{aligned}
&\frac{x^2}{4A} \textrm{ if } x\in]0,2AT[,\\
&Tx-AT^2\geq \frac{Tx}{2} \textrm{ if } x\geq 2AT,
\end{aligned}\right.$$
which ends the proof.
\hfill{\rule{2mm}{2mm}\vskip3mm \par}
\end{rem}

\begin{lem}\label{lemsuperdiff1.2} Let $X$ be a real-valued random variable defined
on some probability space $(\Omega, \cal F, \bb P)$, with $\bb E X\leq 0$ and
$ \bb E[e^{\abs{X}}]\leq K$ for some $K>0$.
Then for all $t\in ]0,1[$,
\begin{equation}
\bb E[e^{tX}]\leq \exp\left(\frac{Kt^2}{1-t}\right).
\end{equation}
Consequently,
\begin{equation}
\bb E[e^{tX}]\leq \exp\left(2Kt^2\right) \textrm{ for every }t\in \left]0,\frac{1}{2}\right].
\end{equation}
\end{lem}

\textbf{Proof.} Let $t\in ]0,1[$. Since $\bb E X\leq 0$, we have
$$\begin{aligned}
\bb E[e^{tX} ]
&=\sum_{k=0}^\infty t^k\bb E\left[\frac{X^k}{k!}\right]
\leq 1+ \sum_{k=2}^\infty t^k\bb E\left[\frac{X^k}{k!} \right]\\
&\leq 1+\sum_{k=2}^\infty t^k\bb E[e^{\abs{X}} ]
\leq 1+K\frac{t^2}{1-t}\leq \exp\left(\frac{Kt^2}{1-t}\right).
\end{aligned}$$
   \hfill{\rule{2mm}{2mm}\vskip3mm \par}

\begin{lem}\label{lemsuperdiff1.3}
For $K>0$ and $x>0$,
\begin{equation}
\sup_{t\in ]0,1[}\left(tx-\frac{Kt^2}{1-t}\right)=\left(\sqrt{x+K}-\sqrt{K}\right)^2.
\end{equation}
\end{lem}

\textbf{Proof.} Let $l_K(t)=\frac{Kt^2}{1-t}$. We first consider
$l_1(t)=\frac{t^2}{1-t}$ (the case where $K=1$). Let
$h(t)=xt-\frac{t^2}{1-t}$, $t\in]0,1[$. Notice that $h'(t)=0$ if and
only if $x=\frac{t(2-t)}{(1-t)^2}$, that is,
$t=1-\frac{1}{\sqrt{1+x}}$. Therefore
$$l_1^*(x)=h\left(1-\frac{1}{\sqrt{1+x}}\right)=(\sqrt{x+1}-1)^2.$$
In the general case, we have
$$l_K^*(x)=Kl_1^*\left(\frac{x}{K}\right)=\left(\sqrt{x+K}-\sqrt{K}\right)^2.$$
\hfill{\rule{2mm}{2mm}\vskip3mm \par}

\textbf{Proof of Theorem \ref{superdiff1}.}
By Lemma \ref{lemsuperdiff1.2}, we obtain that for every $i$ and for every $t\in ]0,1[$, a.s.
$$\bb E[e^{t X_i}|\cal F_{i-1}]\leq \exp\left(\frac{Kt^2}{1-t}\right).$$
Therefore by Lemmas \ref{lemsuperdiff1.1} and \ref{lemsuperdiff1.3},
we obtain immediately (\ref{growthsuperdiff1}) and (\ref{concsuperdiff1}).

To show (\ref{concsuperdiff1bis}), we notice that
the function $g(x)=\frac{\left(\sqrt{x+K}-\sqrt{K}\right)^2}{x^2}$ is strictly decreasing on $]0,+\infty[$
with $\Lim{x}{+\infty}g(x)=0$ and $\Lim{x}{0}g(x)=\frac{1}{4K}$, whereas the function
$f(x)=\frac{\left(\sqrt{x+K}-\sqrt{K}\right)^2}{x}$ is strictly increasing on $]0,+\infty[$,
with $\Lim{x}{+\infty}f(x)=1$ and $\Lim{x}{0}f(x)=0$.

Therefore for every $x\in ]0,K]$, $\left(\sqrt{x+K}-\sqrt{K}\right)^2\geq x^2g(K)=\frac{x^2}{K(1+\sqrt{2})^2}$,
and for every $x>K$, $\left(\sqrt{x+K}-\sqrt{K}\right)^2\geq xf(K)=\frac{x}{(1+\sqrt{2})^2}$,
which ends the proof of (\ref{concsuperdiff1bis}).

Conversely, suppose that $(X_k)$ are iid, and that
$P\left[\frac{S_n}{n}>x\right]\leq e^{-ncx}$ for some $n>1$, $c>0$
and all $x\geq x_1>0$ large enough. Let $\delta \in ]0,c[$.
Then for all $x>0$,
$$\left(P[X_1>x]\right)^n=P[X_i>x \textrm{ for all } 1\leq i\leq n]
\leq P\left[\frac{S_n}{n}>x\right]\leq e^{-ncx},$$
so that $P[X_1^+>x]=P[X_1>x]\leq e^{-cx}$, and
$$\bb E[e^{\delta {X_1^+}}]=1+\int_{0}^{+\infty}P[X_1^+>x]\delta e^{\delta x}dx
\leq 1+\int_{0}^{x_1}\delta e^{\delta x}dx+\int_{x_1}^{+\infty}\delta e^{-(c-\delta) x}dx
=e^{\delta x_1}+\frac{\delta}{c-\delta}e^{-(c-\delta)x_1}.$$
\hfill{\rule{2mm}{2mm}\vskip3mm \par}

\begin{rem}
Notice that by Lemma \ref{lemsuperdiff1.2}, $\forall t\in]0,\frac{1}{2}]$, a.s.
$$\bb E[e^{t X_i}|\cal F_{i-1}]\leq \exp\left(2Kt^2\right).$$
Therefore by Lemma \ref{lempetrov}, we obtain immediately,
\begin{equation}
P\left[\frac{S_n}{n}>x\right]\leq
\left\{\begin{aligned}
&\exp\left(-\ds\frac{nx^2}{8K}\right)\textrm{ if } x\in ]0,2K],\\
&\exp\left(-\ds\frac{nx}{4}\right)\textrm{ if } x>2K.
\end{aligned}\right.
\end{equation}
But (\ref{concsuperdiff1}) of Theorem \ref{superdiff1} gives more precise information.
\end{rem}

\textbf{Proof of Corollary \ref{cor1superdiff1}.}
For $\eps \in ]0,1[$, let $x_0>0$  and $x_1>0$ be such that $g(x_0)=\frac{1}{4K(1+\eps)}$
and $f(x_1)=\frac{1}{1+\eps}$, where $g$ and $f$ are as in the proof of Theorem
\ref{superdiff1}.

If $x\in ]0,x_0]$, then $\left(\sqrt{x+K}-\sqrt{K}\right)^2\geq x^2g(x_0)$, hence
$P\left[\frac{S_n}{n}>x\right] \leq
\exp\left(-\frac{nx^2}{4K(1+\eps)}\right)$.

If $x\in [x_1, +\infty[$, then $\left(\sqrt{x+K}-\sqrt{K}\right)^2\geq xf(x_1)$,
hence $P\left[\frac{S_n}{n}>x\right] \leq
\exp\left(-\frac{nx}{1+\eps}\right)$.

If $x\in [x_0,x_1]$, then $\left(\sqrt{x+K}-\sqrt{K}\right)^2\geq xf(x_0)=xx_0g(x_0)=\frac{xx_0}{4K(1+\eps)}$.
We set $K_1=\frac{4K(1+\eps)}{x_0}$, so that
$P\left[\frac{S_n}{n}>x\right] \leq \exp\left(-\frac{nx}{K_1}\right)$.
\hfill{\rule{2mm}{2mm}\vskip3mm \par}

If we impose an exponential moment condition to $X_i^2$ instead of
$X_i$, we get the following Hoeffding type inequality.

\begin{thm}\label{superdiff2}
Let $(X_i)_{1\leq i\leq n}$ be a sequence of supermartingale
differences adapted to  $(\cal F_i)$. If there exist some constants $R>0$ and $K>0$
such that for all $i$,
\begin{equation}
\bb E[e^{R X_i^2}|\cal F_{i-1}]\leq K \ \ \ a.s.,
\end{equation}
then there exists a constant $c>0$  depending only on $R$ and $K$ such that:
\begin{equation}\label{superdiff2-a}
\bb E[e^{tS_n}]\leq e^{nct^2} \textrm{ for all } t>0,
\end{equation}
and
\begin{equation}\label {superdiff2-b}
 P\left[\frac{S_n}{n}>x\right]\leq e^{-\frac{nx^2}{4c}}
\textrm{ for all } x>0.
\end{equation}
Conversely, if $(X_i)$ are iid and if $(\ref{superdiff2-a})$ or
$(\ref{superdiff2-b})$ holds for some $n\geq 1$ and $c>0$, then for each $R\in \left]0,\frac{1}{4c}\right[$,
$$\bb E[e^{R X_1^{+2}}]\leq K,  \;\;  \mbox{ where }  X_1^+ = \max (X_1, 0)
\textrm{ and }K=1+\frac{R}{\frac{1}{4c}-R}.$$
\end{thm}

Its proof will be based on the following Lemma.
\begin{lem}\label{lemsuperdiff2}
Let $X$ be a random variable defined on a probability space
$(\Omega, \cal F, \bb P)$. If for some constants $K$ and $R>0$,  $\bb
E[e^{R X^2}]\leq K,$ then for all $ t>0$,
\begin{equation}
\bb E[e^{t|X|}]\leq
1+\frac{K\sqrt{\pi}}{2\sqrt{R}}t\exp\left(\frac{t^2}{4R}\right).
\end{equation}
If additionally $\bb E[X]\leq 0$, then there exists $a>0$ depending only on $K$ and $R$ such that for all $t>0$,
\begin{equation}
\bb E[e^{tX}]\leq
\exp\left(\frac{at^2}{2}\right).
\end{equation}
\end{lem}

\textbf{Proof.} By hypothesis $P[\abs{X}>x] \leq
e^{-R x^2}\bb E[e^{R X^2}] \leq K e^{-R x^2}$.   Hence
for all $t>0$,
\begin{eqnarray*}
\bb E[e^{t\abs{X}}] =\int_0^{+\infty}
         P [ e^{t\abs{X}}>x ]dx
= \int_{-\infty}^{+\infty}
         P [\abs{X}>u ]  d(e^{tu})  \nonumber\\
= 1 + t \int_{0}^{+\infty}
         P [\abs{X}>u ]  e^{tu} du
 \leq 1 +  K t \int_{0}^{+\infty}e^{-R u^2}e^{tu}du   \nonumber \\
 \leq 1+\frac{K\sqrt{\pi}}{2\sqrt{R}} t \exp\left(\frac{t^2}{4R}\right).
 \end{eqnarray*}
 Let $c>\frac{1}{4R}$. Then there exists $t_1>0$ such that
 \begin{equation}\label{lemsuperdiff2-1}
 \forall t\geq t_1,\ \ \bb E[e^{t\abs{X}}]
  \leq  \exp\left(ct^2\right).
  \end{equation}
On the other hand,
$$\bb E[e^{R \abs{X}}]\leq
\bb E[e^R; \abs{X}\leq 1]+\bb E[e^{RX^2}; \abs{X}> 1]
\leq e^R+K,$$
so by Lemma \ref{lemsuperdiff1.2}, when $\bb E[X]\leq 0$, we have
\begin{equation}\label{lemsuperdiff2-2}
\bb E[e^{t X} ]\leq \exp\left(\frac{2K_1}{R^2}t^2\right)\ \ \
\forall t\in\left]0,\frac{R}{2}\right],
  \end{equation}
where $K_1=e^R+K$.
From (\ref{lemsuperdiff2-1}) and (\ref{lemsuperdiff2-2}) we deduce that
there exists $a>0$ depending only on $K$ and $R$, such that
\begin{equation*}
 \forall t\geq 0,\ \ \bb E[e^{tX}]\leq
\exp\left(\frac{at^2}{2}\right).
\end{equation*}
\hfill{\rule{2mm}{2mm}\vskip3mm \par}

\textbf{Proof of Theorem \ref{superdiff2}.}
Write $\bb E_{i-1}[.]=\bb E[. |\cal F_{i-1}]$.
By Lemma \ref{lemsuperdiff2} there exists $a=a(R,K)>0$ such that
$$\bb E_{i-1}[e^{t X_i}]\leq \exp\left(\frac{at^2}{2}\right)
\ \ \ \forall t>0.$$
So by Lemmas \ref{lemsuperdiff1.1} and \ref{lempetrov}, we get (\ref{superdiff2-a}) and
(\ref{superdiff2-b}).

Conversely, suppose that $(X_i)$ are iid and that
$(\ref{superdiff2-b})$ holds for some $n\geq 1$ and $c>0$
(notice that $(\ref{superdiff2-a})$ implies $(\ref{superdiff2-b})$).
Let $R \in ]0,\frac{1}{4c}[$.
Then $\forall x>0$,
$$\left(P[X_1>x]\right)^n=P[X_i>x \textrm{ for all } 1\leq i\leq n]
\leq P\left[\frac{S_n}{n}>x\right]\leq e^{-\frac{nx^2}{4c}},$$
so that $P[X_1>x]\leq e^{-\frac{x^2}{4c}}$,
and
$$\bb E[e^{R {X_1^{+2}}}]=1+\int_{0}^{+\infty}P[X_1^+>x]2xR e^{R x^2}dx
=1+\int_{0}^{+\infty}P[X_1>x]2xR e^{R x^2}dx\leq K,$$
where $K=1+\int_{0}^{+\infty}2xR e^{-(\frac{1}{4c}-R) x^2}dx=1+\frac{R}{\frac{1}{4c}-R}$.
\hfill{\rule{2mm}{2mm}\vskip3mm \par}

\section{Exponential bounds of $\bb P (S_n>nx)$ for large values of $x$ }

\setcounter{equation}{0}

Notice that in the exponential inequality $P(S_n\geq nx) \leq
e^{-n c(x)}$  of the preceding section, for large  $x$, we can take
 $c(x) = cx$ or $cx^2$ according to an exponential moment condition
on $X$ or on $X^2$, respectively. In this section we shall see that
this property remains true for $c(x) = cx^Q$ with any $Q\geq 1$.

\begin{thm}\label{suiteadaptee}
Let $(X_i)_{1\leq i\leq n}$ be any adapted sequence with respect to
a filtration $(\cal F_i)_{1\leq i\leq n}$. Assume that there exist some
constants $Q>1$, $R>0$ and $K>0$ such that for all $i\in \interventier{1}{n}$,
\begin{equation}
\bb E[e^{R\abs{X_i}^Q}|\cal F_{i-1}]\leq K\ \ \ a.s..
\end{equation}
Let $\rho>1$ and $\tau>0$ be such that
\begin{equation}
\ds\frac{1}{Q}+\ds\frac{1}{\rho}=1 \textrm{ and }(\rho\tau)^{\frac{1}{\rho}}
(QR)^{\frac{1}{Q}}=1.
\end{equation}
Then for any $\tau_1>\tau$, there exists
$t_1>0$ depending only on $K, Q, R$ and $\tau_1$, such that:
\begin{equation}\label{laplacemdg-ajoute}
\bb E[e^{t\abs{S_n}}]\leq \exp\left(n\tau_1
t^\rho\right)\textrm{ for all } t\geq t_1,
\end{equation}
\begin{equation}\label{tailmdg-ajoute}
P\left[\frac{\abs{S_n}}{n}>x\right]\leq  \exp\left(-nR_1
x^Q\right)\textrm{ for all } x\geq x_1 : = \rho\tau_1t_1^{\rho-1},
\end{equation}
where $R_1>0$ is such that $(\rho\tau_1)^{\frac{1}{\rho}}
(QR_1)^{\frac{1}{Q}}=1.$

Conversely, if $(X_i)$ are iid and if $(\ref{tailmdg-ajoute})$
holds for some $n\geq 1$, $R_1>0$, $Q>1$ and $x_1>0$, then for all $R\in ]0,R_1[$,
\begin{equation}\label{converselyadapted}
\bb E[e^{R \abs{X_1}^{Q}}]\leq 2K,  \;\;  \mbox{ where }
K=e^{Rx_1^Q}+\frac{R}{R_1-R}e^{-(R_1-R)x_1^Q}.
\end{equation}
\end{thm}

When $(X_i)$ are supermartingale differences, we can complete Theorem \ref{suiteadaptee}
with an information for small values of $x>0$ and $t>0$, as shown in the following theorem.
The conclusion follows from Theorem \ref{superdiff1} for small values of $x,t>0$, and from
Theorem \ref{suiteadaptee} for large values of $x,t>0$. The proof of (\ref{converselysupermart})
will be seen in the proof of (\ref{converselyadapted}).
Notice that for large values of $x,t>0$, the conclusion of Theorem \ref{superdiffQ} is
sharper than that of Theorem \ref{superdiff1}.

\begin{thm}\label{superdiffQ}
Under the hypothesis of Theorem \ref{suiteadaptee}, if
moreover $(X_i)_{1\leq i\leq n}$ is a sequence of supermartingale
differences  adapted to the filtration $(\cal F_i)$,
then for any $\tau_1>\tau$, there exist $t_1>0$, $x_1>0$, and $A,B>0$,
depending only on $K, Q, R$, and $\tau_1$, such that:
\begin{equation}\label{laplacemdg}
\bb E[e^{tS_n}]\leq \left\{\begin{aligned}
&\exp\left(n\tau_1 t^\rho\right)\textrm{ if } t\geq t_1,\\
&\exp\left(nAt^2\right)\textrm{ if } 0\leq t\leq t_1,
\end{aligned}\right.
\end{equation}
and \begin{equation}\label{tailmdg}
P\left[\frac{S_n}{n}>x\right]\leq \left\{\begin{aligned}
&\exp\left(-nR_1 x^Q\right)\textrm{ if } x\geq x_1,\\
&\exp\left(-nBx^2\right)\textrm{ if } 0\leq x\leq x_1.\\
\end{aligned}\right.
\end{equation}
Conversely, if $(X_i)$ are iid and if the first inequality in $(\ref{tailmdg})$
holds for some $n\geq 1$, $R_1>0$, $Q>1$ and $x_1>0$, then for all $R\in ]0,R_1[$,
\begin{equation}\label{converselysupermart}
\bb E[e^{R X_1^{+Q}}]\leq K,  \;\;  \mbox{ where }  X_1^+ = \max
(X_1, 0) \textrm{ and }
K=e^{Rx_1^Q}+\frac{R}{R_1-R}e^{-(R_1-R)x_1^Q}.
\end{equation}
\end{thm}

Before proving the theorems, we first give, for a positive random
variable $X$,  relations among
 the growth rate of the Laplace transform
$\bb E[e^{tX}]$ (as $t\rightarrow\infty$),  the decay rate of the tail
probability  $P [X>x]$ (as $x\rightarrow \infty$), and the
exponential moments of the form $\bb E[e^{RX^Q}]$ $(Q>1)$.

\begin{lem}\label{lemfond1} (Relation between $\bb E[e^{tX}]$ and $P [X>x]$)

Let $X$ be a
positive real random variable. Let $Q$, $\rho$, $\tau$, and $R\in
]0,+\infty[$ be such that $1<Q<+\infty$ and
$$\ds\frac{1}{Q}+\ds\frac{1}{\rho}=1,\ \ \ (\rho\tau)^{\frac{1}{\rho}} (QR)^{\frac{1}{Q}}=1.$$
Let $K>0$ be a constant. Consider the following assertions:

\medskip
\begin{list}{}
\item (1) $\forall t>0,\ \bb E[e^{tX}]\leq Ke^{\tau t^ \rho}$;
\item (2) $\forall x>0,\ P[X>x]\leq Ke^{-Rx^Q} $;
\item (3) For $a=K\left(\frac{2}{R}\right)^{\frac{1}{Q-1}}$ and all $t>0$,
$\bb E[e^{tX}]\leq 1+K+at^\rho e^{\tau t^ \rho}$.
\end{list}
Then we have the following implications: $(1)\imp (2)\imp (3)$.
\end{lem}

Lemma \ref{lemfond1} is closely related to the following Legendre
duality between the functions $t\mapsto \tau t^\rho$ and $x\mapsto
Rx^Q$.

\begin{lem}\label{lemopt}
Let $\rho>1$, $\tau>0$ and $t_0\geq 0$. Then $\forall x\geq \rho\tau t_0^{\rho-1}$,
$$\sup_{t\geq t_0}\left(tx-\tau t^\rho\right)=Rx^Q,\textrm{ where }\ds\frac{1}{Q}+\ds\frac{1}{\rho}=1,
(\rho\tau)^{\frac{1}{\rho}} (QR)^{\frac{1}{Q}}=1.$$
\end{lem}

\textbf{Proof.}
The fonction $h(t)=tx-\tau t^\rho$ attains its supremum on $]0,+\infty[$ for
$t^\star=(\frac{x}{\tau \rho})^{\frac{1}{\rho -1}}$, and the supremum is
$h(t^*)=Rx^Q$. As $t^\star\geq t_0$ if and only if $x\geq  \rho\tau t_0^{\rho-1}$, we get the result.
\hfill{\rule{2mm}{2mm}\vskip3mm \par}

\textbf{Proof of Lemma \ref{lemfond1}.}
We first prove the implication
$(1)\imp (2)$. If $\bb E[e^{tX}]\leq Ke^{\tau t^ \rho}$ then for every $x>0$ and $t>0$,
$$P[X>x]=P[e^{tX}>e^{tx}]\leq e^{-tx} \bb E[e^{tX}]
\leq Ke^{\tau t^\rho-tx}.$$
Therefore by Lemma \ref{lemopt}, $P[X>x] \leq Ke^{-Rx^Q}$.

We then prove the implication $(2)\imp (3)$. If (2) holds, then for every $t>0$,
$$
\bb E[e^{tX}]=
 1+t \int_0^{+\infty} P[X>x]e^{tx}dx
 \leq  1+tK \int_0^{+\infty} e^{-Rx^Q+tx}dx.
$$
We choose $x_1=(\frac{2t}{R})^{\frac{1}{Q -1}}$ so that
$-Rx^Q+tx\leq -xt$ for $x\geq x_1$;
by Lemma \ref{lemopt} (with $t_0=0$), $-Rx^Q+tx\leq \tau t^\rho$ for any $x>0$. Therefore
$$ \int_0^{+\infty} e^{-Rx^Q+tx}dx\leq  \int_0^{x_1} e^{\tau t^\rho}dx
+ \int_{x_1}^{+\infty} e^{-xt}dx
\leq x_1e^{\tau t^\rho}+\frac{1}{t},$$
hence for $a=K(\frac{2}{R})^{\frac{1}{Q-1}}$ and $t>0$,
$$\bb E[e^{tX}]\leq
1+K+at^{\rho}e^{\tau t^\rho}.$$
\hfill{\rule{2mm}{2mm}\vskip3mm \par}

\begin{lem}\label{lemfond2}
Let $X$ be a positive real random variable.
Let $Q\in [1,+\infty[$, and $K$, $R\in ]0,+\infty[$.
Consider the following assertions:
\begin{list}{}{}
\item (1) $\bb E[e^{RX^Q}]\leq K$;
\item (2) $\forall x>0$, $P[X>x]\leq Ke^{-Rx^Q}$;
\item (3) For any $R_1\in ]0,R[$, $\bb E[e^{R_1X^Q}]\leq \frac{R+R_1(K-1)}{R-R_1}$.
\end{list}
Then we have the following implications: $(1)\imp (2)\imp (3)$.
\end{lem}

\textbf{Proof of Lemma \ref{lemfond2}.}
The implication $(1)\imp (2)$ is easy: if $\bb E[e^{RX^Q}]\leq K$, then
$P[X>x]=P[e^{RX^Q}>e^{Rx^Q}]\leq Ke^{-Rx^Q}$.
Let us now prove the implication $(2)\imp (3)$.
If $P[X>x]\leq Ke^{-Rx^Q}$, then for any $R_1\in ]0,R[$,
$$
\begin{aligned}
\bb E[e^{R_1X^Q}]=
\int_0^{+\infty} P[e^{R_1X^Q}>x] dx
&= 1+R_1 Q \int_0^{+\infty} P[X>u]e^{R_1 u^Q}u^{Q-1}du\\
&\leq  1+KR_1 Q \int_0^{+\infty} e^{(R_1-R)u^Q}u^{Q-1}du
=\ds\frac{R+R_1(K-1)}{R-R_1}.
\end{aligned}$$
\hfill{\rule{2mm}{2mm}\vskip3mm \par}

\begin{rem}
Let $Q$, $\rho\in ]0,+\infty[$ be such that $1<Q<+\infty$ and
$\frac{1}{Q}+\frac{1}{\rho}=1$.
As a consequence of Lemma \ref{lemfond1}, we can easily see
that writing
$$\begin{aligned}
\tau&=\inf\{a>0:\ \bb E[e^{rX}]=O(\exp(ar^\rho))\},\\
R&=\sup\{a>0:\ P[X>x]=O(\exp(-ax^Q))\},
\end{aligned}$$
we have
$$(\rho\tau)^{\frac{1}{\rho}} (QR)^{\frac{1}{Q}}=1.$$
This was proved in a different way by Liu in \cite{Liu1996}.
It unifies Theorems 6.1, 7.1, 7.2, 7.3, 8.1, 9.1 and 9.2 of Ramachandran (\cite{Ramachandran1962}),
and was first conjectured by Harris (\cite{Harris1948}) in the context of branching processes.
\end{rem}

\textbf{Proof of Theorem \ref{suiteadaptee}.} By Lemmas
\ref{lemfond2} and \ref{lemfond1}, we see that for
$a=K\left(\frac{2}{R}\right)^{\frac{1}{Q-1}}$,
$$\bb E[e^{t\abs{X_i}}|\cal F_{i-1}]\leq 1+K+at^\rho e^{\tau t^\rho}\ \ \ \forall t>0.$$
Let $\tau_1>\tau$. Then  there exists $t_1>0$ sufficiently
large such that $\forall t\geq t_1$, $\bb E[e^{t\abs{X_i}}|\cal F_{i-1}]\leq e^{\tau_1t^\rho}$.
Applying Lemmas \ref{lemsuperdiff1.1} and \ref{lemopt},  we obtain that
$$\bb E[e^{t\abs{S_n}}]
\leq \bb E[e^{t(\abs{X_1}+\cdots+\abs{X_n})}]
\leq\exp\left(n\tau_1 t^\rho\right)\textrm{   if } t\geq t_1,$$
$$P\left[\frac{\abs{S_n}}{n}>x\right]
\leq P\left[\frac{\abs{X_1}+\cdots+\abs{X_n}}{n}>x\right]
\leq e^{-nR_1x^Q} \textrm{   if }x\geq x_1=\rho\tau_1 t_1^{\rho-1}.$$
Conversely, suppose that $(X_k)$ are iid, and that
$P\left[\frac{\abs{S_n}}{n}>x\right]\leq \exp\left(-nR_1
x^Q\right)\textrm{ for all } x\geq x_1$.
Let $R\in ]0,R_1[$.
Then for all $x\geq x_1$,
$$\left(P[X_1>x]\right)^n=P[X_i>x \textrm{ for all } 1\leq i\leq n]
\leq P\left[\frac{S_n}{n}>x\right]\leq P\left[\frac{\abs{S_n}}{n}>x\right]\leq \exp\left(-nR_1 x^Q\right),$$
so that $X_1^+=\max(0,X_1)$ satisfies $P[X_1^+>x]=P[X_1>x]\leq e^{-R_1x^Q}$,
and
\begin{multline*}
\bb E[e^{R (X_1^{+})^Q}]=1+\int_{0}^{+\infty}P[X_1^+>x]R Qx^{Q-1}e^{R x^Q}dx\\
\leq 1+\int_{0}^{x_1}R Qx^{Q-1}e^{R x^Q}dx+\int_{x_1}^{+\infty}R Qx^{Q-1}e^{-(R_1-R) x^Q}dx
=e^{Rx_1^Q}+\frac{R}{R_1-R}e^{-(R_1-R)x_1^Q}.
\end{multline*}
By considering $(-S_n)$ instead of $(S_n)$, we see that the same result holds for $X_1^-=\max(0,-X_1)$:
$$\bb E[e^{R {(X_1^{-})}^Q}]\leq K:=e^{Rx_1^Q}+\frac{R}{R_1-R}e^{-(R_1-R)x_1^Q}.$$
Therefore $$\bb E[e^{R \abs{X_1}^{Q}}]\leq \bb E[e^{R (X_1^{+})^
Q}]+\bb E[e^{R (X_1^{-})^Q}]\leq 2K.$$

\hfill{\rule{2mm}{2mm}\vskip3mm \par}

\textbf{Proof of Theorem \ref{superdiffQ}.}
By Theorem \ref{suiteadaptee}, there exists $t_1>\frac{R}{2}$ such that
$$\bb E[e^{tS_n}]\leq \bb E[e^{t\abs{S_n}}]\leq \exp\left(n\tau_1
t^\rho\right)\textrm{ for all } t\geq t_1,$$
and $$P\left[\frac{S_n}{n}>x\right]\leq P\left[\frac{\abs{S_n}}{n}>x\right]\leq  \exp\left(-nR_1
x^Q\right)\textrm{ for all } x\geq x_1 : = \rho\tau_1t_1^{\rho-1}.$$
On the other hand, notice that $\bb E[e^{R\abs{X_i}}|\cal F_{i-1}]\leq K_1:=e^R+K$,
so that by Theorem \ref{superdiff1},
$$\bb E[e^{tS_n}]\leq \exp\left(\frac{2nK_1t^2}{R^2}\right)
\textrm{ for all }t\in \left]0,\frac{R}{2}\right].$$
If $t\in\left]\frac{R}{2},t_1\right]$, then
$$\bb E[e^{tS_n}]\leq \bb E[e^{t_1\abs{S_n}}]\leq \exp\left(n\tau_1 t_1^\rho\right)
\leq e^{n\frac{4\tau_1 t_1^\rho}{R^2}t^2}.$$
Set $A=\max(\frac{2K_1}{R^2},\frac{4\tau_1 t_1^\rho}{R^2})$.
Then $$\bb E[e^{tS_n}]\leq e^{nAt^2} \ \forall t\in ]0,t_1].$$
Again by Theorem \ref{superdiff1}, we can choose $B>0$ small enough such that
$$P\left[\frac{S_n}{n}>x\right]\leq
\exp\left(-nBx^2\right)\textrm{ if }x\in ]0,x_1].$$
\hfill{\rule{2mm}{2mm}\vskip3mm \par}

\section{Extension to the case $\bb E[e^{\abs{X_i}}|\cal F_{i-1}]\leq K_i$}
\setcounter{equation}{0}

 The following theorems are immediate
generalizations of Theorems \ref{superdiff1}, \ref{superdiff2},
\ref{suiteadaptee} and \ref{superdiffQ}. The proofs of the first two
theorems remain the same; the proof of the third needs a short
argument for the concerned constants to be independent of $n$. The
first theorem is an extension of Bernstein's inequality.

   \begin{thm} \label{superdiff1Ki}
Let $(X_i)_{1\leq i\leq n}$ be a finite sequence of supermartingale
differences. If for some constants $K_i>0$ and all $i\in\interventier{1}{n}$, a.s.
\begin{equation}
\bb E[e^{\abs{X_i}}|\cal F_{i-1}]\leq K_i,
\end{equation}
then for each $K\geq \frac{K_1+\cdots +K_n}{n}$,
\begin{equation}
 \bb E[e^{tS_n}]\leq \exp\left({\frac{nKt^2}{1-t}}\right)
\textrm{ for all }t\in ]0,1[,
\end{equation}
and \begin{equation}
P\left[\frac{S_n}{n}>x\right]\leq\exp\left(-n\left(\sqrt{x+K}-\sqrt{K}\right)^2\right) \textrm{ for all } x>0.
\end{equation}
Consequently,
 \begin{equation}
P\left[\frac{S_n}{n}>x\right]\leq \left\{\begin{aligned}
&\exp\left(-\ds\frac{nx^2}{K(1+\sqrt{2})^2}\right)\textrm{ if } x\in ]0,K],\\
&\exp\left(-\ds\frac{nx}{(1+\sqrt{2})^2}\right)\textrm{ if } x>K.
\end{aligned}\right.
\end{equation}
\end{thm}

The second theorem is an extension of Hoeffding's inequality.

\begin{thm}\label{superdiff2Ki}
Let $(X_i)_{1\leq i\leq n}$ be a sequence of supermartingale
differences adapted to  $(\cal F_i)$. If there exist some constants $K_i>0$
and $R>0$ such that for all $i\in\interventier{1}{n}$, a.s.
\begin{equation}
\bb E[e^{R X_i^2}|\cal F_{i-1}]\leq K_i,
\end{equation}
then for each $K\geq \frac{K_1+\cdots +K_n}{n}$, there exists
a constant $c>0$  depending only on $R$ and $K$ such that:
\begin{equation}
\bb E[e^{tS_n}]\leq e^{nct^2} \textrm{ for all } t>0,
\end{equation}
and
\begin{equation}
 P\left[\frac{S_n}{n}>x\right]\leq e^{-\frac{nx^2}{4c}}
\textrm{ for all } x>0.
\end{equation}
\end{thm}

The third theorem shows a close relation between $P[\abs{X_i}>x]$
and $P\left[\frac{\abs{S_n}}{n}>x\right]$ for large values of $x>0$. Notice that this
result is valid for any adapted sequence.

\begin{thm}\label{suiteadapteeKi}
Let $(X_i)_{1\leq i\leq n}$ be any adapted sequence with respect to
a filtration $(\cal F_i)_{1\leq i\leq n}$. Assume that there exist some
constants $Q>1$, $R>0$ and $K_i>0$ such that for all $i\in \interventier{1}{n}$,
\begin{equation}
\bb E[e^{R\abs{X_i}^Q}|\cal F_{i-1}]\leq K_i\ \ \ a.s..
\end{equation}
Let $\rho>1$ and $\tau>0$ be such that
\begin{equation}
\ds\frac{1}{Q}+\ds\frac{1}{\rho}=1 \textrm{ and }(\rho\tau)^{\frac{1}{\rho}}
(QR)^{\frac{1}{Q}}=1.
\end{equation}
Let $K\geq \frac{K_1+\cdots +K_n}{n}$.
Then for any $\tau_1>\tau$, there exists
$t_1>0$ depending only on $K, Q, R$ and $\tau_1$, such that:
\begin{equation}\label{laplacemdg-ajouteKi}
\bb E[e^{t\abs{S_n}}]\leq \exp\left(n\tau_1
t^\rho\right)\textrm{ for all } t\geq t_1,
\end{equation}
and \begin{equation}\label{tailmdg-ajouteKi}
P\left[\frac{\abs{S_n}}{n}>x\right]\leq  \exp\left(-nR_1
x^Q\right)\textrm{ for all } x\geq x_1 : = \rho\tau_1t_1^{\rho-1},
\end{equation}
where $R_1$ is such that $(\rho\tau_1)^{\frac{1}{\rho}}
(QR_1)^{\frac{1}{Q}}=1.$
\end{thm}

\textbf{Proof.} By Lemmas
\ref{lemfond2} and \ref{lemfond1}, we see that for
$a=\left(\frac{2}{R}\right)^{\frac{1}{Q-1}}$,
$$\bb E[e^{t\abs{X_i}}|\cal F_{i-1}]\leq 1+K_i(1+at^\rho e^{\tau t^\rho})\ \ \ \forall t>0.$$
By Lemma \ref{lemsuperdiff1.1},
$$\bb E[e^{t\abs{S_n}}]
\leq \bb E[e^{t(\abs{X_1}+\cdots+\abs{X_n})}]
\leq \prod_{i=1}^ {n}\left(1+K_i(1+at^\rho e^{\tau t^\rho})\right).$$
It is easy to see that $1+K_i(1+at^\rho e^{\tau t^\rho})\leq e^{K_i}(1+at^\rho e^{\tau t^\rho})$, so we have
$$\bb E[e^{t\abs{S_n}}]\leq \left(e^K (1+at^\rho e^{\tau t^\rho})\right)^n.$$
Let $\tau_1>\tau$. Then  there exists $t_1>0$ sufficiently
large such that $\forall t\geq t_1$, $e^K (1+at^\rho e^{\tau t^\rho})\leq e^{\tau_1t^\rho}$,
which gives (\ref{laplacemdg-ajouteKi}). As
$$P\left[\frac{\abs{S_n}}{n}>x\right]=P[e^{t\abs{S_n}}> e^{tnx}]
\leq e^{-ntx}\bb E[e^{t\abs{S_n}}]\leq \exp\left(-n(tx-\tau_1t^\rho)\right),$$
we deduce (\ref{tailmdg-ajouteKi}) from Lemma \ref{lemopt}.
\hfill{\rule{2mm}{2mm}\vskip3mm \par}

As in section 3, when $(X_i)$ are supermartingale differences, using Theorem \ref{superdiff1Ki}
we can complete Theorem \ref{suiteadapteeKi}
with an information for small values of $x>0$ and $t>0$, as shown in the following theorem.
For large values of $x,t>0$, it gives inequalities sharper than those of Theorem \ref{superdiff1Ki}.

\begin{thm}\label{superdiffQKi}
Under the same hypothesis as in Theorem \ref{suiteadapteeKi}, if
moreover $(X_i)_{1\leq i\leq n}$ is a sequence of supermartingale
differences  adapted to the filtration $(\cal F_i)$,
then for any $\tau_1>\tau$, there exist $t_1>0$, $x_1>0$, and $A,B>0$,
depending only on $K, Q, R$ and $\tau_1$, such that:
\begin{equation}\label{laplacemdgKi}
\bb E[e^{tS_n}]\leq \left\{\begin{aligned}
&\exp\left(n\tau_1 t^\rho\right)\textrm{ if } t\geq t_1,\\
&\exp\left(nAt^2\right)\textrm{ if } 0\leq t\leq t_1,
\end{aligned}\right.
\end{equation}
and \begin{equation}\label{tailmdgKi}
P\left[\frac{S_n}{n}>x\right]\leq \left\{\begin{aligned}
&\exp\left(-nR_1 x^Q\right)\textrm{ if } x\geq x_1,\\
&\exp\left(-nBx^2\right)\textrm{ if } 0\leq x\leq x_1.\\
\end{aligned}\right.
\end{equation}
\end{thm}

\section{Rate of convergence with probability 1 and in $L^p$}
\setcounter{equation}{0}

   \begin{thm} \label{superdiffps}
Let $(X_i)_{1\leq i\leq n}$ be a sequence of supermartingale
differences. If for some constants $K_i>0$
and for all $i\in\interventier{1}{n}$,
\begin{equation}\label{Ki}
\bb E[e^{\abs{X_i}}|\cal F_{i-1}]\leq K_i \ \ \ a.s.,
\end{equation}
then writing $K=\limsup_{n\to +\infty}\frac{K_1+\cdots +K_n}{n}$ and
$S_n^+ = \max(0, S_n)$, we have:
\begin{equation}\label{limsupps}
\limsup_{n\to +\infty} \frac{S_n^+}{\sqrt{n\ln n}}\leq 2\sqrt K \ \ \ a.s.,
\end{equation}
and for every $p>0$,
\begin{equation}\label{growthLp}
\limsup_{n\to +\infty}\ \ n^{\frac{p}{2}}\bb E\left[\left(\frac{S^+_n}{n}\right)^p\right]\leq
 p2^{p-1}K^{\frac{p}{2}}\Gamma\left(\frac{p}{2}\right).
\end{equation}
 \end{thm}

\textbf{Proof.} For the proof of (\ref{limsupps}), by Borel-
Cantelli's Lemma, it suffices to show that for every $a>2\sqrt K$,
$$\Sum{n=0}{+\infty}P\left[\frac{S_n^+}{\sqrt{n\ln n}}> a\right]<+\infty.$$
Let us fix $a>2\sqrt K$. Let $\eps>0$ be such that $a>2\sqrt{K+\eps}$
and let $n_1>0$ be such that for every $n\geq n_1$,
$\frac{K_1+\cdots +K_n}{n}<K+\eps$.
Then we deduce from Theorem \ref{superdiff1Ki} that for every $n\geq n_1$,
$$\begin{aligned}
P\left[\frac{S_n^+}{\sqrt{n\ln n}}> a\right]=P\left[\frac{S_n}{n}>\sqrt{\frac{\ln n}{n}}a\right]
\leq \exp\left(-n\left(\sqrt{x_n+K+\eps}-\sqrt{K+\eps}\right)^2\right)
\end{aligned},$$
with $x_n=\sqrt{\frac{\ln n}{n}}a$.
When $n$ tends to $\infty$,
$n\left(\sqrt{x_n+K+\eps}-\sqrt{K+\eps}\right)^2\sim \frac{a^2\ln n}{4(K+\eps)}$.
As $a^2>4(K+\eps)$, it follows that
$$\Sum{n=0}{+\infty}P\left[\frac{S_n^+}{\sqrt{n\ln n}}> a\right]<+\infty.$$
This ends the proof of (\ref{limsupps}).

We now come to the proof of (\ref{growthLp}).
Let $n_1>0$ be as in the proof of (\ref{limsupps}).
We deduce from Theorem \ref{superdiff1Ki} that for every $n\geq n_1$,
$$\begin{aligned}
\bb E\left[\left(\frac{S^+_n}{n}\right)^p\right]&=p\int_0^\infty P\left[\frac{S_n}{n}>x\right]x^{p-1}dx
\leq p\int_0^\infty \exp\left(-n\left(\sqrt{x+K+\eps}-\sqrt{K+\eps}\right)^2\right)x^{p-1}dx.
\end{aligned}$$
Set $y=n\left(\sqrt{x+K+\eps}-\sqrt{K+\eps}\right)^2$. Then
$\sqrt{\frac{y}{n}}=\sqrt{x+K+\eps}-\sqrt{K+\eps}$,
$x=\sqrt{\frac{y}{n}}\left(\sqrt{\frac{y}{n}}+2\sqrt {K+\eps}\right)$, $dx=\frac{\sqrt{\frac{y}{n}}+\sqrt{K+\eps}}{\sqrt{n}\sqrt{y}}dy$, so that
\begin{equation}\label{}
\bb E\left[\left(\frac{S^+_n}{n}\right)^p\right]
\leq \frac{c_n(p)}{n^{\frac{p}{2}}},
\end{equation}
where $$c_n(p)=p\int_0^\infty e^{-y}y^{\frac{p}{2}-1}
\left(\sqrt{\frac{y}{n}}+2\sqrt {K+\eps}\right)^{p-1}\left(\sqrt{\frac{y}{n}}+\sqrt{K+\eps}\right)dy$$
satisfies
$$\lim_{n\to \infty}c_n(p)=p\int_0^\infty e^{-y}y^{\frac{p}{2}-1}
\left(2\sqrt {K+\eps}\right)^{p-1}\left(\sqrt{K+\eps}\right)dy=
p2^{p-1}(K+\eps)^{\frac{p}{2}}\Gamma\left(\frac{p}{2}\right).$$
\hfill{\rule{2mm}{2mm}\vskip3mm \par}

In the case of a sequence of martingale differences, replacing $S_n^+$ by $\abs{S_n}$
in the proof above, we obtain immediately:

 \begin{cor} \label{diffps}
Let $(X_i)_{1\leq i\leq n}$ be a sequence of martingale
differences. If for some constants $K_i>0$
and for all $i\in\interventier{1}{n}$,
\begin{equation}\label{condmartdiffKi}
 \bb E[e^{\abs{X_i}}|\cal F_{i-1}]\leq K_i \ \ \ a.s.,
\end{equation}
then for $K=\limsup_{n\to +\infty}\frac{K_1+\cdots +K_n}{n}$,
\begin{equation}\label{limsuppsmart}
\limsup_{n\to +\infty} \frac{\abs{S_n}}{\sqrt{n\ln n}}\leq 2\sqrt K  \ \ \ a.s.,
\end{equation}
and for every $p>0$,
\begin{equation}\label{growthLpmart}
\limsup_{n\to +\infty}\ \ n^{\frac{p}{2}}\bb E\left[\left(\frac{\abs{S_n}}{n}\right)^p\right]\leq
 p2^{p}K^{\frac{p}{2}}\Gamma\left(\frac{p}{2}\right).
\end{equation}
 \end{cor}

\begin{rem}
The exponential moment condition (\ref{condmartdiffKi}) can certainly be relaxed for a result
of the form $\bb E\left[\left(\frac{\abs{S_n}}{n}\right)^p\right]=O\left(n^{-\frac{p}{2}}\right)$.
For example, as shown in \cite{LesigneVolny}, p.150, by Burkholder's inequality, we can obtain the
following result:
if $p\geq 2$ and $\bb E[\abs{X_i}^p]\leq K$ for some $K>0$ and all $i\in\interventier{1}{n}$, then
\begin{equation}
\bb E[\abs{S_n}^p]\leq n^{\frac{p}{2}}(18pq^{1/2})^pK,
\end{equation}
where $\frac{1}{p}+\frac{1}{q}=1$.
\end{rem}

\section{Free energy of directed polymers: concentration inequalities}
\setcounter{equation}{0}

We now consider the model of a directed polymer in a random environment, already described in
the introduction. For convenience, let us recall it briefly as follows.
Let $\omega=(\omega_n)_{n\in\bb N}$ be the simple random walk on the $d$-dimensional integer lattice $\bb Z^d$ starting at $0$, defined on a probability space $(\Omega,\cal F,P)$.
Let $\eta=(\eta(n,x))_{(n,x)\in\bb N\times \bb Z^d}$ be a sequence of real valued, non constant and i.i.d. random variables defined on another probability space $(E,\cal E, \bb Q)$.
The path $\omega$ represents the directed polymer and $\eta$ the random environment.
For any $n>0$, define the random polymer measure $\mu_n$ on the path space $(\Omega,\cal F)$ by
\begin{equation}
\mu_n=\ds\frac{1}{Z_n(\beta)}\exp(\beta H_n(\omega))P(d\omega),
\end{equation}
where $\beta\in\bb R$ is the inverse temperature,
\begin{equation}
H_n(\omega)=\sum_{j=1}^n \eta(j,\omega_j),\textrm{ and } Z_n(\beta)=P[\exp(\beta H_n(\omega))].
\end{equation}
Let $\lambda(\beta)=\ln \bb Q[e^{\beta\eta(0,0)}]$ be the logarithmic moment generating function of
$\eta(0,0)$. We fix $\beta >0$ (otherwise we consider $-\eta$), and assume only $\lambda(\pm \beta)<\infty$,
which is equivalent to $\bb Q [e^{\beta|\eta(0,0)|}]<\infty$.
We are interested in the asymptotic behaviour of the normalized partition function
\begin{equation}
W_n(\beta)=\frac{Z_n(\beta)}{\bb Q[Z_n(\beta)]}=P[\exp(\beta H_n-n\lambda(\beta))],
\end{equation}
and the free energy $\frac{1}{n}\ln W_n(\beta)$.
For simplicity, we shall write $W_n$ for $W_n(\beta)$, $Z_n$ for $Z_n(\beta)$, and $\eta$ for $\eta(0,0)$.
We use the same letter $\eta$ to denote the environment sequence $(\eta(n,x))_{(n,x)\in\bb N\times \bb Z^d}$
and the random variable $\eta(0,0)$; there will be no confusion according to the context.
In this section, we shall prove exponential concentration inequalities for the free energies
$\frac{\ln W_n}{n}$, and convergence results of the centered energies
$\frac{\ln W_n}{n}-\frac{\bb Q[\ln W_n]}{n}$: cf. Theorems \ref{concentration1},
\ref{estimation}, \ref{concentrationQ}, and their corollaries.

\begin{thm}\label{concentration1}
Assume that $\bb Q[e^{\beta\abs{\eta}}]<+\infty$, and
set $K=2\exp\left(\lambda(-\beta))+\lambda(\beta)\right)$. Then for all $n\geq 1$,
\begin{equation}\label{growth}
\bb Q[e^{\pm t(\ln W_n-\bb Q[\ln W_n])}]\leq \exp\left(\frac{nKt^2}{1-t}\right) \textrm{ for all }t\in ]0,1[,
\end{equation}
and \begin{equation}\label{speed}
\bb Q\left[\pm\frac{1}{n}(\ln W_n-\bb Q[\ln W_n])>x\right]
\leq \exp\left(-n\left(\sqrt{x+K}-\sqrt{K}\right)^2\right) \textrm{ for all } x>0.
\end{equation}
Consequently, $\forall n\geq 1$,
\begin{equation}\label{speedbis}
 \bb Q\left[\pm\frac{1}{n}(\ln W_n-\bb Q[\ln W_n])>x\right]\leq \left\{\begin{aligned}
&\exp\left(-\ds\frac{nx^2}{K(1+\sqrt{2})^2}\right)\textrm{ if } x\in ]0,K],\\
&\exp\left(-\ds\frac{nx}{(1+\sqrt{2})^2}\right)\textrm{ if } x>K.
\end{aligned}\right.
\end{equation}
\end{thm}

\begin{cor}\label{cor1concentration1}
Under the conditions of Theorem \ref{concentration1}, $\forall\eps >0$, there exist $0<x_0<x_1$ and $K_1>0$ depending only on $K$ and $\eps$, such that:
\begin{equation}\label{corspeed1}
\bb Q\left[\pm\frac{1}{n}(\ln W_n-\bb Q[\ln W_n])>x\right] \leq
\left\{\begin{aligned}
&\exp\left(-\ds\frac{nx^2}{4K(1+\eps)}\right)\textrm{ if } x\in ]0,x_0[,\\
&\exp\left(-\ds\frac{nx}{K_1}\right)\textrm{ if } x\in [x_0,x_1],\\
&\exp\left(-\ds\frac{nx}{1+\eps}\right)\textrm{ if } x\in ]x_1,+\infty[.
\end{aligned}\right.
\end{equation}
\end{cor}

\begin{rem}
Using Lesigne and Volny's martingale inequality (\ref{LV}),
Comets, Shiga and Yoshida (2003, \cite{CometsShigaYoshida}) proved that
if $\bb Q[e^{\beta\abs{\eta}}]<+\infty$ for all $\beta>0$, then
$\forall x>0$, there exists $n_0\in\bb N^*$ such that for any $n\geq n_0$,
\begin{equation}
\bb Q\left[\abs{\frac{1}{n}\ln W_n-\frac{1}{n}\bb Q[\ln W_n]}>x\right]\leq
\exp\left(-\ds\frac{n^{\frac{1}{3}}x^{\frac{2}{3}}}{4}\right).
\end{equation}
Our result is sharper as $n^{1/3}$ is replaced by $n$.
Another advantage is that our conclusion holds for all $n$, not only for $n$ large enough;
thanks to this advantage, we can use our inequalities to study the convergence rate for the a.s. and $L^p$ convergence:
cf. Theorem \ref{estimation}.
The third advantage is that we assume $\bb Q[e^{\beta\abs{\eta}}]<+\infty$
only  for the fixed $\beta$, not for all $\beta>0$.
The first two advantages are due to the application of our exponential martingale inequality
(Theorem \ref{superdiff1}); the third one comes from a direct estimation of the conditional
exponential moment (Lemma \ref{CSY})
by use of convex inequalities, without using Lemma 3.1 of \cite{CometsShigaYoshida}.
\end{rem}

For the proof, as in \cite{CometsShigaYoshida},
we write $\ln W_n-\bb Q[\ln W_n]$ as a sum of $(\cal E_j)_{1\leq j\leq n}$ martingale differences:
$$\ln W_n-\bb Q[\ln W_n]=\ds\sum_{j=1}^nV_{n,j}, \textrm{ with } V_{n,j}=\bb Q_j[\ln W_n]-\bb Q_{j-1}[\ln W_n],$$
where $\bb Q_j$ denotes the conditional expectation with respect to $\bb Q$ given $\cal E_j$,
$\cal E_j=\sigma[\eta(i,x): 1\leq i\leq j, x\in \bb Z^d]$.

\begin{lem}\label{CSY} We have
\begin{equation}\label{inegCSY}
\bb Q_{j-1}\left [ \exp(tV_{n,j})\right ]\leq \exp\left(L(t)\right)\textrm{ for every } t\in\bb R,
\end{equation}
where
\begin{equation}
L(t)=\left\{\begin{aligned}
&\lambda(t\beta)+\lambda(-t\beta) \textrm{ if }\abs{t}>1,\\
&\lambda(-\abs{t}\beta)+\abs{t}\lambda(\beta) \textrm{ if }\abs{t}\leq 1.
\end{aligned}\right.
\end{equation}
Consequently,
\begin{equation}\label{CSYbis}
\bb Q_{j-1}\left [ \exp(\abs{V_{n,j}})\right]\leq
K:= 2\exp\left(\lambda(\beta)+\lambda(-\beta)\right).
\end{equation}
\end{lem}

\textbf{Proof.}
We fix $t\in\bb R^*$ and assume $L(t)<\infty$ (otherwise there is nothing to prove).
Set
$$e_{n,j}=\exp\left (\sum_{1\leq k\leq n,k\neq j} (\beta\eta(k,\omega_k)-\lambda(\beta))\right ),
\ \ \ W_{n,j}=P[e_{n,j}].$$
Since $\bb Q_{j-1}[\ln W_{n,j}]=\bb Q_j[\ln W_{n,j}],$ we have
\begin{equation}\label{W}
V_{n,j}=\bb Q_j\left [\ln \frac{W_n}{W_{n,j}}\right ]-
\bb Q_{j-1}\left [\ln \frac{W_n}{W_{n,j}}\right ].
\end{equation}
For $j\in\bb N$ and $x\in\bb Z^d$, define
$$\overline{\eta}_x=\overline{\eta}(j,x)=\exp(\beta \eta(j,x)-\lambda(\beta)),\ \ \alpha_x=\frac{P[e_{n,j};\omega_j=x]}{W_{n,j}}.$$
(Throughout the paper, for a measure $\mu$, a function $f$, and a set $A$, we use the notation
$\mu[f;A]=\int f\mathbf 1_Ad\mu$, where $\mathbf 1_A$ is the indicator function of $A$).
Then $$\sum_{x\in \bb Z^d}\alpha_x=1 \textrm{ and }\ds\frac{W_n}{W_{n,j}}=
\sum_{x\in \bb Z^d}\alpha_x\overline{\eta}_x.$$
By (\ref{W}),
\begin{equation*}
\bb Q_{j-1}\left [ \exp(t V_{n,j})\right ]=
\exp\left(-t \bb Q_{j-1} \left [\ln \frac{W_n}{W_{n,j}}\right ]\right)
\bb Q_{j-1}\left [\exp\left(t \bb Q_j \left [\ln \frac{W_n}{W_{n,j}}\right ]\right)\right].
\end{equation*}
Since the function $x\mapsto e^{tx}$ is convex,
using Jensen's inequality and the fact that $\cal E_{j-1}\subset\cal E_j$, we get:
\begin{equation}\label{CSY1}
\bb Q_{j-1}\left [ \exp(t V_{n,j})\right ]\leq
\bb Q_{j-1}\left[\left(\frac{W_n}{W_{n,j}}\right)^{-t}\right]
\bb Q_{j-1} \left[\left( \frac{W_n}{W_{n,j}}\right)^t\right].
\end{equation}
If $t<0$ or $t\geq 1$ then the function $x\mapsto x^{t}$ is convex, therefore by
 Jensen's inequality we have
$$\left(\frac{W_n}{W_{n,j}}\right)^{t}=
\left(\sum_{x\in \bb Z^d}\alpha_x\overline{\eta}_x\right)^{t}\leq
\sum_{x\in \bb Z^d}\alpha_x\left(\overline{\eta}_x\right)^{t}.$$
We consider
the $\sigma$-algebra $\cal E_{n,j}=\sigma[\eta(k,x); 1\leq k\leq n, k\neq j,  x\in \bb Z^d]$.
Then $\cal E_{j-1}\subset \cal E_{n,j}$, the $\alpha_x$
are $\cal E_{n,j}$-measurable,
and the $\overline{\eta}_x$ are independent of $\cal E_{n,j}$, so that
\begin{equation*}
\bb Q_{j-1} [\alpha_x\left(\overline{\eta}_x\right)^{t}]
=\bb Q_{j-1} [\bb Q[\alpha_x\left(\overline{\eta}_x\right)^{t} |\cal E_{n,j}]]
=\bb Q_{j-1} [\alpha_x\bb Q[\left(\overline{\eta}_x\right)^{t}]]
=\exp\left(\lambda(t\beta)-t\lambda(\beta)\right)\bb Q_{j-1} [\alpha_x].
\end{equation*}
Hence for $t<0$ or $t\geq 1$,
\begin{equation}\label{CSY2}
\bb Q_{j-1}\left[\left(\frac{W_n}{W_{n,j}}\right)^{t}\right]
\leq
\exp\left(\lambda(t\beta)-t\lambda(\beta)\right).
\end{equation}
It is easily seen that the equality holds for $t=1$: $\bb Q_{j-1}\left[\frac{W_n}{W_{n,j}}\right]=1$.
Again by Jensen's inequality, we have, for $t\in ]0,1]$,
\begin{equation}\label{CSY3}
\bb Q_{j-1}\left[\left(\frac{W_n}{W_{n,j}}\right)^{t}\right]
\leq \left(\bb Q_{j-1}\left[\frac{W_n}{W_{n,j}}\right]\right)^{t}=1.
\end{equation}
The inequality (\ref{inegCSY}) is then just a combination of (\ref{CSY1}), (\ref{CSY2}), and (\ref{CSY3}).
In particular,
$$\bb Q_{j-1}\left [ \exp(\pm V_{n,j})\right]\leq \exp\left(\lambda(\beta)+\lambda(-\beta)\right),
\textrm{ so that }
\bb Q_{j-1}\left [ \exp(\abs{V_{n,j}})\right]\leq
K:= 2\exp\left(\lambda(\beta)+\lambda(-\beta)\right).$$
\hfill{\rule{2mm}{2mm}\vskip3mm \par}

\textbf{Proof of Theorem  \ref{concentration1}.}
From Lemma \ref{CSY} and Theorem \ref{superdiff1}, we deduce:
\begin{equation}\label{growth+}
\bb Q[e^{t(\ln W_n-\bb Q[\ln W_n])}]\leq \exp\left(\frac{nKt^2}{1-t}\right) \textrm{ for all }t\in ]0,1[,
\end{equation}
and
\begin{equation}\label{speed+}
 \bb Q\left[\frac{1}{n}\left(\ln W_n-\bb Q[\ln W_n]\right)>x\right]
\leq \exp\left(-n\left(\sqrt{x+K}-\sqrt{K}\right)^2\right) \textrm{ for all } x>0.
\end{equation}
Applying Theorem \ref{superdiff1} to the sequence $(-V_{n,j})$, we find that
\begin{equation}\label{growth-}
\bb Q[e^{-t(\ln W_n-\bb Q[\ln W_n])}]\leq \exp\left(\frac{nKt^2}{1-t}\right) \textrm{ for every }t\in ]0,1[,
\end{equation}
and
\begin{equation}\label{speed-}
 \bb Q\left[-\frac{1}{n}\left(\ln W_n-\bb Q[\ln W_n]\right)>x\right]
\leq \exp\left(-n\left(\sqrt{x+K}-\sqrt{K}\right)^2\right) \textrm{ for all } x>0.
\end{equation}
The inequalities (\ref{growth+}) and (\ref{growth-}) give
(\ref{growth}), (\ref{speed+}) and (\ref{speed-}) give
(\ref{speed}). \hfill{\rule{2mm}{2mm}\vskip3mm \par}

\textbf{Proof of Corollary \ref{cor1concentration1}.}
The proof is the same as the proof of Corollary \ref{cor1superdiff1}.
\hfill{\rule{2mm}{2mm}\vskip3mm \par}

\begin{thm}\label{estimation}
Assume that $\bb Q[e^{\beta\abs{\eta}}]<+\infty$, and
set $K=2\exp\left(\lambda(-\beta))+\lambda(\beta)\right)$.
Then
\begin{equation}\label{convergenceenergy}
\frac{1}{n}\ln W_n-\frac{1}{n}\bb Q[\ln W_n]\to 0 \ \ a.s.\textrm{ and in } L^p,
\end{equation}
with
\begin{equation}\label{estimationps}
\limsup_{n\to +\infty} \ \sqrt{\frac{n}{\ln n}}\abs{\frac{\ln W_n}{n}-\frac{\bb Q[\ln W_n]}{n}}\leq 2\sqrt K\ \ a.s.,
\end{equation}
and for every $p>0$,
\begin{equation}\label{estimationLp}
\limsup_{n\to +\infty}\ \ n^{\frac{p}{2}}\bb Q\left[\abs{\frac{\ln W_n-\bb Q[\ln W_n]}{n}}^p\right]\leq
p2^{p}K^{\frac{p}{2}}\Gamma\left(\frac{p}{2}\right).
\end{equation}
\end{thm}

\textbf{Proof.}
Recall that with the notations of the proof of Theorem \ref{concentration1},
we have
$$\bb Q_{j-1}\left [ \exp(\abs{V_{n,j}})\right ]\leq K.$$
Then the inequalities (\ref{estimationps}) and (\ref{estimationLp})
are consequences of the inequalities  (\ref{limsuppsmart}) and (\ref{growthLpmart}) of Corollary \ref{diffps}.
\hfill{\rule{2mm}{2mm}\vskip3mm \par}

\begin{thm}\label{concentrationQ}
Assume that
$K_0:=\bb Q[e^{R\abs{\eta}^Q}]<+\infty$ for some $Q>1$
and $R>0$.
Let $\rho>1$ and $\tau>0$ be determined by
\begin{equation}
\ds\frac{1}{Q}+\ds\frac{1}{\rho}=1 \textrm{ and }(\rho\tau)^{\frac{1}{\rho}}
(QR)^{\frac{1}{Q}}=1.
\end{equation}
Then for each $\tau_1>\tau$, there exist constants $t_0, A, B>0$, depending only on
$\beta$, $K_0$, $Q$, $R$, $\tau$ and $\tau_1$, such that, for all $n\geq 1$,
\begin{equation}\label{laplacemdQ}
\bb Q[e^{\pm t(\ln W_n-\bb Q[\ln W_n])}]\leq \left\{\begin{aligned}
&\exp\left(2n\tau_1\beta^\rho t^\rho\right)\textrm{ if } t> \frac{t_0}{\beta},\\
&\exp (nAt^2)\textrm{ if } 0< t\leq \frac{t_0}{\beta},
\end{aligned}\right.
\end{equation}
and \begin{equation}\label{tailmdQ}
\bb Q\left[\pm\frac{1}{n}(\ln W_n-\bb Q[\ln W_n])>x\right]\leq \left\{\begin{aligned}
&\exp\left(-nR_1 x^Q\right)\textrm{ if } x> 2\rho\beta\tau_1 t_0^{\rho-1},\\
&\exp\left(-nBx^2\right)\textrm{ if } 0< x\leq 2\rho\beta\tau_1t_0^{\rho-1},\\
\end{aligned}\right.
\end{equation}
where $R_1>0$ is such that $\beta(2\rho\tau_1)^{\frac{1}{\rho}} (QR_1)^{\frac{1}{Q}}=1$.
\end{thm}
If we are not interested in the values of constants, then we have

\begin{cor}\label{corollaireQ}
Under the conditions of Theorem \ref{concentrationQ},
there exist constants $c_1,c_2>0$, depending only on $\beta$, $K_0$, $Q$ and $R$, such that:
\begin{equation}
\bb Q[e^{\pm t(\ln W_n-\bb Q[\ln W_n])}]\leq \left\{\begin{aligned}
&\exp\left(nc_1 t^\rho\right)\textrm{ if } t> 1,\\
&\exp (nc_1t^2)\textrm{ if } 0< t\leq 1,
\end{aligned}\right.
\end{equation}
and \begin{equation}
\bb Q\left[\pm\frac{1}{n}(\ln W_n-\bb Q[\ln W_n])>x\right]\leq \left\{\begin{aligned}
&\exp\left(-nc_2 x^Q\right)\textrm{ if } x> 1,\\
&\exp\left(-nc_2x^2\right)\textrm{ if } 0< x\leq 1.\\
\end{aligned}\right.
\end{equation}
In particular, if $K_0:=\bb Q[e^{R\abs{\eta}^2}]<+\infty$ for some
$R>0$, then for some constants $c_1,c_2>0$ depending only on $\beta$, $K_0$ and $R$,
\begin{equation}\label{conc2}
\bb Q[e^{\pm t(\ln W_n-\bb Q[\ln W_n])}]\leq
\exp\left(nc_1 t^2\right)\textrm{ for all } t\in\bb R,
\end{equation}
and \begin{equation}
\bb Q\left[\pm\frac{1}{n}(\ln W_n-\bb Q[\ln W_n])>x\right]\leq
\exp\left(-nc_2 x^2\right)\textrm{ for all } x> 0.
\end{equation}
\end{cor}

\begin{rem}
If the environment is bounded or gaussian,
the inequality (\ref{conc2})
was proved in \cite{CometsVargas}, Corollary 2.5,
as a corollary of a general concentration  result.
\end{rem}

\textbf{Proof of Theorem \ref{concentrationQ}.}
Let $\tau_1>\tau$.
By Lemmas \ref{lemfond2} and \ref{lemfond1}, writing $a=K_0(\frac{2}{R})^{\frac{1}{Q-1}}$, we have
$$\bb Q[e^{t\abs{\eta}}]\leq 1+K_0+at^\rho e^{\tau t \rho}\leq e^{\tau_1 t^\rho}\ \forall t\geq t_0,$$
for some $t_0=t_0(K_0,\rho,\tau,\tau_1)>1$. Hence $\lambda(\pm t)\leq \tau_1 t^\rho$ $\forall t\geq t_0,$
so by Lemma \ref{CSY},
\begin{equation}
\bb Q_{j-1}\left [ \exp(\pm tV_{n,j})\right ]\leq \exp\left(L(\pm  t)\right)\leq
\exp\left(2\tau_1\beta^\rho t^\rho\right)\textrm{ for all } t\geq \frac{t_0}{\beta}.
\end{equation}
We apply Lemma \ref{lemsuperdiff1.1} with $I=]\frac{t_0}{\beta},+\infty[$, and with the aid of Lemma \ref{lemopt},
we conclude that
\begin{equation}\label{laplacemdQ0}
\bb Q[e^{\pm t(\ln W_n-\bb Q[\ln W_n])}]\leq \exp\left(2n\tau_1 \beta^\rho t^\rho\right)
\textrm{ if } t> \frac{t_0}{\beta},
\end{equation}
and
\begin{equation}\label{tailmdQ0}
\bb Q\left[\pm \frac{1}{n}(\ln W_n-\bb Q[\ln W_n])>x\right]\leq
\exp\left(-nR_1 x^Q\right)\textrm{ if } x> 2\rho \tau_1 \beta t_0^{\rho-1}.
\end{equation}
Clearly, the condition $\bb Q[e^{R|\eta|^Q}] < \infty$ implies $\bb
Q[e^{\beta|\eta|}] < \infty$. Let
 $K=2\exp\left(\lambda(\beta)+\lambda(-\beta)\right)$.
 By Theorem \ref{concentration1},
(\ref{growth}), and Corollary \ref{cor1concentration1},
(\ref{corspeed1}), $\forall\eps>0$,
\begin{equation}\label{laplacemdQ1}
\bb Q[e^{\pm t(\ln W_n-\bb Q[\ln W_n])}]\leq \exp\left(\frac{nKt^2}{1-t}\right)
\leq \exp\left(\frac{nKt^2}{1-\eps}\right) \textrm{ if }0<t\leq\eps,
\end{equation}
and \begin{equation}\label{tailmdQ1}
\bb Q\left[\pm\frac{1}{n}(\ln W_n-\bb Q[\ln W_n])>x\right]
\leq \exp\left(-\frac{nx^2}{4K(1+\eps)}\right) \textrm{ if } 0<x<\delta(K,\eps),
\end{equation}
for some $\delta(K,\eps)$ small enough. In the following, we take
$\eps=\frac{1}{2}$ and $\delta=\delta(K,\frac{1}{2})$. If
$\frac{1}{2}\leq t\leq \frac{t_0}{\beta}$, then
\begin{equation}\label{laplacemdQ2}
\bb Q[e^{\pm t(\ln W_n-\bb Q[\ln W_n])}]\leq
\bb Q\left[e^{\frac{t_0}{\beta}\abs{\ln W_n-\bb Q[\ln W_n]}}\right]\leq
2\exp\left(2n\tau_1 t_0^\rho\right)
\leq \exp\left(n(4\ln 2+8\tau_1t_0^\rho)t^2\right).
\end{equation}
Combining (\ref{laplacemdQ0}), (\ref{laplacemdQ1}) and (\ref{laplacemdQ2}) gives
(\ref{laplacemdQ}), with $A=\max(4K,4\ln 2+8\tau_1t_0^\rho)$.

If $\delta\leq x\leq x_0:=2\rho \tau_1 \beta t_0^{\rho-1}$, then by (\ref{speed}),
\begin{equation}\label{tailmdQ2}
\bb Q\left[\pm \frac{1}{n}(\ln W_n-\bb Q[\ln W_n])>x\right]\leq
\exp\left(-n(\sqrt{\delta +K}-\sqrt{K})^2\right)\leq
\exp\left(-n\frac{(\sqrt{\delta +K}-\sqrt{K})^2x^2}{x_0^2}\right).
\end{equation}
Combining (\ref{tailmdQ0}), (\ref{tailmdQ1}) and (\ref{tailmdQ2}) gives (\ref{tailmdQ}),
with $B=\min\left(\frac{1}{6K}, \frac{(\sqrt{\delta +K}-\sqrt{K})^2}{x_0^2}\right)$.
\hfill{\rule{2mm}{2mm}\vskip3mm \par}

\section{Free energy of directed polymers: convergence rates}
\setcounter{equation}{0}

It is well known that the sequence $\bb Q[\ln W_n(\beta)]$ is
superadditive, hence the limit
\begin{equation}\label{superadd}
p_{-}(\beta)=\Lim n \infty \frac{1}{n}\bb Q[\ln (W_n(\beta))]=
\sup_n \frac{1}{n}\bb Q[\ln (W_n(\beta))]\in ]-\infty,0]
\end{equation}
exists{\footnote {In the literature, $p(\beta)$ is often used to
denote the limit of the un-normalized free energy: $p(\beta)=\Lim n
\infty \frac{1}{n}\bb Q[\ln (Z_n(\beta))].$ We use the symbol
$p_{-}(\beta)$ to indicate that $p_{-}(\beta)\leq 0$. Of course
$p_{-}(\beta)=p(\beta)-\lambda(\beta)$.}}.
 As an
immediate consequence of (\ref{superadd}) and
(\ref{convergenceenergy}), we have:

\begin{lem}\label{}
Assume that $\bb Q[e^{\beta\abs{\eta}}]<+\infty$. Then
\begin{equation}
p_{-}(\beta)=\Lim n \infty \frac{1}{n}\ln (W_n(\beta)) \in \,
[\beta\bb Q[\eta]-\lambda(\beta), \; 0] \quad Q\mbox{-a.s. and in }
L^p, \;\; \forall p\geq 1.
\end{equation}
\end{lem}

The inequality $p_{-}(\beta)\leq 0$ was already indicated in (\ref{superadd});
it follows from the fact that
$$\bb Q[\ln (W_n(\beta))]\leq \ln\bb Q[W_n(\beta)]=0.$$
The inequality $p_{-}(\beta)\geq \beta\bb Q[\eta]-\lambda(\beta)$
also comes directly from the definition, as
$$\bb Q[\ln (W_n)]\geq \bb Q P[\beta H_n-n\lambda(\beta)]
=P\bb Q[\beta H_n-n\lambda(\beta)]=n(\beta\bb
Q[\eta]-\lambda(\beta)).$$ The a.s. convergence was proved in
\cite{CometsShigaYoshida}, under the stronger condition that $\bb
Q[e^{\beta\abs{\eta}}]<+\infty$ for all $\beta >0$; actually their
proof is valid under the condition that $\bb
Q[e^{3\beta\abs{\eta}}]<+\infty$. We shall give an estimation of the
rate of convergence, for each of the convergences in probability,
a.s., and in $L^p$ ($p\geq 1$): cf. Theorems \ref{energy} and
\ref{estimationdetermin}.

We first consider the rate of convergence in probability. Recall
that the condition $\bb Q[e^{\beta\abs{\eta}}]<+\infty$ is
equivalent to $\lambda(\pm \beta) <\infty$.

\begin{thm}\label{energy}
If $K:=2\exp\left(\lambda(-\beta))+\lambda(\beta)\right) <\infty$,
then $\forall \delta\in]0,1[$, $\forall x>0$, there exists
$n_0=n_0(\delta,x)>0$ such that $\forall n\geq n_0$,
\begin{equation}\label{speedenergy}
\bb Q\left[\abs{\frac{1}{n}\ln W_n-p_{-}(\beta)}>x\right] \leq
2\exp\left(-n\left(\sqrt{(1-\delta)x+K}-\sqrt{K}\right)^2\right).
\end{equation}
Consequently,
\begin{equation}\label{speedenergybis}
 \bb Q\left[\abs{\frac{1}{n}\ln W_n-p_{-}(\beta)}>x\right]\leq \left\{\begin{aligned}
&2\exp\left(-\ds\frac{n(1-\delta)^2x^2}{K(1+\sqrt{2})^2}\right)\textrm{ if } x(1-\delta)\leq K,\\
&2\exp\left(-\ds\frac{n(1-\delta)x}{(1+\sqrt{2})^2}\right)\textrm{ if } x(1-\delta)>K.
\end{aligned}\right.
\end{equation}
In particular (take $\delta=\frac{1}{2}$), $\forall x\in ]0,2K]$, there exists $n_0=n_0(x)>0$ such that $\forall n\geq n_0$,
\begin{equation}
\bb Q\left[\abs{\frac{1}{n}\ln W_n-p_{-}(\beta)}>x\right]\leq
2\exp\left(-\ds\frac{nx^2}{4K(1+\sqrt 2)^2}\right).
\end{equation}
\end{thm}

\textbf{Proof.}
Let $\delta\in]0,1[$, and $x>0$.
Let $n_0=n_0(\delta,x)$ be large enough such that for any $n\geq n_0$,
$$0\leq p_{-}(\beta)-\frac{1}{n}\bb Q[\ln (W_n(\beta))] < \delta x.$$
Then $\forall n\geq n_0$,
$$\bb Q\left[\abs{\frac{1}{n}\ln W_n-p_{-}(\beta)}>x\right]\leq
\bb Q\left[\abs{\frac{1}{n}\ln W_n-\frac{1}{n}\bb Q[\ln W_n]}>(1-\delta)x \right].$$
Therefore the conclusion follows from Theorem \ref{concentration1}.
\hfill{\rule{2mm}{2mm}\vskip3mm \par}

We next consider the rate of convergence in mean. To this end, we first introduce some notations.
We note $P^x$ the law of the simple random walk on
$\bb Z^d$ starting at $x$, and $L_m=\{x\in\bb Z^d,\ P(\omega_m=x)>0\}$.
In addition to the partition function $W_n$, we define the partition function starting from $x$:
\begin{equation}
W_{n}(x)=W_{n}(x;\eta)=P^x\left[\exp\left(\beta \sum_{j=1}^{n} \eta(j,\omega_j)-n\lambda(\beta\right)\right],
\end{equation}
and the point to point partition function
\begin{equation}\label{partitionfunction}
W_{n}(x,y)=W_{n}(x,y;\eta)=P^x\left[\exp\left(\beta \sum_{j=1}^{n} \eta(j,\omega_j)-n\lambda(\beta)\right)\textbf{1}_{\omega_{n}=y}\right].
\end{equation}
Let $\tau_n$ be the time shift of ordre $n$ on the environment:
$$(\tau_n\eta)(k,x)=\eta(k+n,x) \ \ \ (x\in\bb Z^d,k\geq 1).$$
Then we have
\begin{equation}\label{markov}
W_{n+k}=\sum_{x\in L_n}W_n(0,x;\eta)W_k(x;\tau_n\eta).
\end{equation}

\begin{lem}\label{speeddeterministe} (Rate of convergence in mean)
If $K:=2\exp\left(\lambda(-\beta))+\lambda(\beta)\right) <\infty$,
then for each $n\in\bb N^*$,
\begin{equation}
0\leq p_{-}(\beta)- \frac{1}{n}\bb Q[\ln (W_n(\beta))]\leq
2\sqrt{K}\sqrt{\frac{d\ln(2n)}{n}}+\frac{d\ln(2n)}{n}.
\end{equation}
\end{lem}

\textbf{Proof.}
We adapt the proof of Proposition 2.4 of \cite{CarmonaHu2004}.
Let $\eps\in ]0,1[$.
Using (\ref{markov}) and  the subadditivity of the function $u\mapsto u^\eps$ ,
we get
$$W_{n+k}^\eps\leq \sum_{x\in L_n}W_n^\eps(0,x;\eta)W_k^\eps (x;\tau_n\eta).$$
Integrating with respect to the environment we have
$$\bb Q[W_{n+k}^\eps]\leq \abs{L_n}\bb Q[W_n^\eps]\bb Q[W_k^\eps]
\leq (2n)^d \bb Q[W_n^\eps]\bb Q[W_k^\eps].$$
Therefore $h_\eps(n):=\ln \bb Q[W_n^\eps]$ $(\geq \eps \bb Q[\ln W_n]$) satisfies
$$h_\eps(n+k)\leq h_\eps(n)+h_\eps(k)+d\ln(2n)\ \ \ \forall n,k\geq 1.$$
Set $h(\eps)=\limsup_{n\to +\infty} \frac{h_\eps(n)}{n}$ (in fact by Hammersley's (1962, \cite{Hammersley1962}) theorem
on sub-additive functions, the limit exists, although we shall not use it).
By the preceding recurrence relation,
we have
$$h_\eps(nm)\leq mh_\eps(n)+(m-1)d\ln(2n),\ \ \ n,m\geq 1.$$
Dividing this inequality by $nm$ and letting $m\to \infty$, we see that
$$h(\eps)\leq \frac{h_\eps(n)+d\ln(2n)}{n},\ \ \forall n\geq 1.$$
By Theorem \ref{concentration1},
$$h_\eps(n)=\ln \bb Q[\exp\left (\eps(\ln W_n-\bb Q[\ln W_n])\right)]
+\eps \bb Q[\ln W_n]\leq \frac{nK\eps^2}{1-\eps}+\eps \bb Q[\ln W_n].$$
As $p_{-}(\beta)\leq \frac{h(\eps)}{\eps}$, it follows that
\begin{equation}\label{deterministe1}
p_{-}(\beta)\leq \frac{K\eps}{1-\eps}+\frac{\bb Q[\ln W_n]}{n}+\frac{d\ln(2n)}{n\eps}
\ \ \ \forall \eps\in ]0,1[.
\end{equation}
Let $g(\eps)=\frac{K\eps}{1-\eps}+\frac{d_n}{\eps}$, where $d_n=\frac{d\ln(2n)}{n}$,
$\eps\in ]0,1[$.
Then $g'(\eps)=\frac{K}{(1-\eps)^2}-\frac{d_n}{\eps^2}=0$ if and only if
$\eps=\frac{\sqrt{d_n}}{\sqrt{K}+\sqrt{d_n}}$.
For $\eps=\frac{\sqrt{d_n}}{\sqrt{K}+\sqrt{d_n}}$, $g(\eps)=2\sqrt{Kd_n}+d_n$; therefore
taking $\eps=\frac{\sqrt{d_n}}{\sqrt{K}+\sqrt{d_n}}$ in (\ref{deterministe1}), we obtain
$$p_{-}(\beta)\leq 2\sqrt{Kd_n}+d_n+\frac{\bb Q[\ln W_n]}{n},$$
that is,
$$p_{-}(\beta)\leq 2\sqrt{K}\sqrt{\frac{d\ln(2n)}{n}}+\frac{d\ln(2n)}{n}+\frac{\bb Q[\ln W_n]}{n}.$$
\hfill{\rule{2mm}{2mm}\vskip3mm \par}

As an immediate consequence of the preceding lemma, we have:
\begin{cor}\label{cordeterm}
 If $K:=2\exp\left(\lambda(-\beta))+\lambda(\beta)\right) <\infty$,
 then
\begin{equation}\label{cordeterm1}
\limsup_{n\to +\infty} \ \sqrt{\frac{n}{\ln n}}\abs{\frac{\bb Q[\ln W_n]}{n}-p_{-}(\beta)}\leq 2\sqrt {Kd}.
\end{equation}
\end{cor}

We finally consider the rate of convergence, with probability 1 and
in $L^p$. As usual, $\N{.}_p$ denotes the $L^p$ norm.

\begin{thm}\label{estimationdetermin}(Rate of convergence, a.s. and in $L^p$)
If $K:=2\exp\left(\lambda(-\beta))+\lambda(\beta)\right) < \infty$,
then
\begin{equation}\label{estimationdeterminps}
\limsup_{n\to +\infty} \ \sqrt{\frac{n}{\ln n}}\abs{\frac{\ln
W_n}{n}-p_{-}(\beta)} \leq 2\sqrt K(1+\sqrt{d})\ \ a.s.,
\end{equation}
and
\begin{equation}\label{estimationdeterminLp}
\limsup_{n\to +\infty}\ \sqrt{\frac{n}{\ln n}}\N{\frac{\ln
W_n}{n}-p_{-}(\beta)}_p \leq 2\sqrt{Kd},  \quad \forall p\geq 1.
\end{equation}

\end{thm}

\textbf{Proof.}
We write
$$\frac{\ln W_n}{n}-p_{-}(\beta)=\left(\frac{\ln W_n}{n}-\frac{\bb Q[\ln W_n]}{n}\right)
+\left(\frac{\bb Q[\ln W_n]}{n}-p_{-}(\beta)\right).$$
Then combining (\ref{estimationps}) of Theorem \ref{estimation} and (\ref{cordeterm1})
of Corollary \ref{cordeterm}, we get (\ref{estimationdeterminps}).
Again by Theorem \ref{estimation}, we know that for every $p\geq 1$,
\begin{equation*}
\N{\frac{\ln W_n-\bb Q[\ln W_n]}{n}}_p=O(n^{-1/2})=o\left(\sqrt{\frac{\ln n}{n}}\right),
\end{equation*}
so that (\ref{estimationdeterminLp}) is a consequence of Corollary
\ref{cordeterm}. \hfill{\rule{2mm}{2mm}\vskip3mm \par}

\begin{rem}
Carmona and Hu have proved in \cite{CarmonaHu2004} that if the environment is gaussian, then for any
$\eps>0$,
$$\abs{\frac{\ln W_n}{n}-p_{-}(\beta)}\leq n^{-(\frac{1}{2}-\eps)} \textrm{ for } n \textrm{ big enough. }$$
Our estimation is sharper since $n^{(\frac{1}{2}-\eps)}$ is replaced by $\sqrt{\frac{\ln n}{n}}$.
\end{rem}

\section{Expression of the free energy by multiplicative cascades}
\setcounter{equation}{0}

In this section we shall prove that the free energy $p_{-}(\beta)$
can be expressed in terms of the free energies of some generalized
multiplicative cascades. The expression is interesting because we
know more information on the free energies of multiplicative
cascades. The model of multiplicative cascades was first introduced
by Mandelbrot (1974, \cite{Man74}); it  has been  well studied in
the literature: see for example Kahane and Peyri\`ere (1976,
\cite{KP76}), Durrett and Liggett (1981,\cite{DL83}), Guivarc'h
(1990, \cite{Gui90}), Franchi (1993, \cite{Fra93} );  for a
generalized version and closely related topics, see Liu (2000,
\cite{Liu2000}).

In \cite{CometsVargas}, Comets and Vargas introduced a generalized
multiplicative cascade (cf. \cite{Liu2000}) $(W_{m,n}^{tree})_{n\geq
1}$ associated to the random vector $(W_m(0,x))_{x\in L_m}$, where
we recall that
\begin{equation}
W_{m}(0,x)=P[\exp(\beta H_{m}(\omega)-m\lambda(\beta));\omega_{m}=x].
\end{equation}
The associated free energy is
\begin{equation}
p_m^{tree}(\beta)=\inf_{\theta\in]0,1]}v_m(\theta),\textrm{ with }
v_m(\theta)=\frac{1}{\theta}\ln\left(\bb Q[\sum_{x\in L_m}W_m(0,x)^\theta]\right).
\end{equation}
 Comets et Vargas proved that
\begin{equation}\label{inegCV}
p_{-}(\beta)\leq \inf_{m\geq 1}\frac{1}{m}p_m^{tree}(\beta)=\limf{m}{+\infty}\frac{1}{m}p_m^{tree}(\beta),
\end{equation}
and that the equality holds if the environment is gaussian or bounded.
Here we prove that the equality holds for general environment.

\begin{thm}\label{cascades}
Assume that $\bb Q[e^{\beta\abs{\eta}}]<+\infty$. Then
\begin{equation}
p_{-}(\beta)= \inf_{m\geq 1}\frac{1}{m}p_m^{tree}(\beta).
\end{equation}
\end{thm}

\textbf{Proof.}
For the sake of completeness, we recall the argument of Comets-Vargas for the inequality
(\ref{inegCV}).
Using the point to point partition functions defined by (\ref{partitionfunction}), we have
\begin{equation}
W_{mn}=\sum_{x_1,\cdots,x_n\in \bb Z^d}W_m(0,x_1)W_m(x_1,x_2;\tau_m\eta)\cdots
W_m(x_{n-1},x_n;\tau_{(n-1)m}\eta).
\end{equation}
Let $\theta\in]0,1[$ and $m\in\bb N^*$.
By the subadditivity of the function $u\mapsto u^\theta$ and Jensen's inequality,
we obtain:
$$\begin{aligned}
\frac{1}{nm}\bb Q[\ln W_{nm}]&=
\frac{1}{m\theta n}\bb Q[\ln W_{nm}^\theta]\\
&\leq \frac{1}{m\theta n}\bb Q\left[\ln \sum_{x_1,\cdots,x_n\in \bb Z^d}W_m^\theta(0,x_1)
 W_m^\theta(x_1,x_2;\tau_m\eta)\cdots
W_m^\theta(x_{n-1},x_n;\tau_{(n-1)m}\eta)\right]\\
&\leq \frac{1}{m\theta n}\ln \bb Q\left[\sum_{x_1,\cdots,x_n\in \bb Z^d}W_m^\theta(0,x_1) W_m^\theta(x_1,x_2;\tau_m\eta)\cdots
W_m^\theta(x_{n-1},x_n;\tau_{(n-1)m}\eta)\right].
\end{aligned}$$
By induction on $n$ it is easy to see that
$$\bb Q\left[\sum_{x_1,\cdots,x_n\in \bb Z^d}W_m^\theta(0,x_1) W_m^\theta(x_1,x_2;\tau_m\eta)\cdots
W_m^\theta(x_{n-1},x_n;\tau_{(n-1)m}\eta)\right]
=\left(\bb Q\left[\sum_{x\in L_m}W_m^\theta(0,x) \right]\right)^n,$$
therefore
$$\frac{1}{nm}\bb Q[\ln W_{nm}]\leq
\frac{1}{m\theta}\ln\left(\bb Q[\sum_{x\in L_m}W_m^\theta(0,x) ]\right).$$
Letting $n\to \infty$ gives
 $$p_{-}(\beta)\leq
\frac{1}{m\theta}\ln\left(\bb Q[\sum_{x\in L_m}W_m^\theta(0,x)
]\right);$$ then taking the infimum over all $\theta\in ]0,1]$ gives
$$p_{-}(\beta)\leq
\frac{1}{m}p_m^{tree}(\beta).$$

Now we prove the reverse inequality.
As $W_m(0,x)\leq W_m$ for every $x$, we have, for $\theta\in]0,1[$,
$$v_m(\theta)\leq\frac{1}{\theta}\ln\left(\bb Q[\abs{L_m}W_m^\theta]\right),$$
where $\abs{L_m}$ is the cardinality of $L_m$. Writing
$$\bb Q[W_m^\theta]=e^{\theta \bb Q[\ln W_m]}
\bb Q[\exp\left (\theta(\ln W_m-\bb Q[\ln W_m])\right)],$$
we get
$$v_m(\theta)
\leq
\frac{1}{\theta}\ln\abs{L_m}+\bb Q[\ln W_m]+\frac{1}{\theta}\ln\left(\bb Q[\exp\left (\theta(\ln W_m-\bb Q[\ln W_m])\right )]\right).$$
Recall that by Theorem \ref{concentration1},
for every $\theta\in ]0,1[$,
$$\bb Q[\exp(\theta (\ln W_m-\bb Q[\ln W_m]))] \leq  e^{\frac{mK\theta^2}{1-\theta}},$$
with $K=2\exp\left(\lambda(-\beta))+\lambda(\beta)\right)$.
Therefore for any $m\geq 1$,
$$\inf_{m\geq 1}\frac{1}{m}p_m^{tree}(\beta)\leq \frac{1}{m}p_m^{tree}(\beta)\leq \frac{1}{m}v_m(\theta)
\leq \frac{1}{m\theta}\ln\abs{L_m}+\frac{1}{m}\bb Q[\ln W_m]+\frac{K\theta}{1-\theta}.$$
Letting $m\to \infty$ and using the fact that $\abs{L_m}\leq (2m)^d$,
we obtain that
$$\inf_{m\geq 1}\frac{1}{m}p_m^{tree}(\beta)\leq p_{-}(\beta)+\frac{K\theta}{1-\theta}.$$
This gives the desired result as $\theta\in]0,1[$ is arbitrary.
\hfill{\rule{2mm}{2mm}\vskip3mm \par}

\bibliographystyle{plain}
\bibliography{biblioprob}

\end{document}